\documentclass[11pt]{article}

\input xy

\xyoption{all}

\usepackage{amssymb,amsbsy,amsthm,amsmath,graphicx,epsfig,times}

\newtheorem{thm}{Theorem}[section]
\newtheorem{lem}[thm]{Lemma}
\newtheorem{prop}[thm]{Proposition}
\newtheorem{cor}[thm]{Corollary}

\newtheorem{dfn}[thm]{Definition}

\theoremstyle{remark}

\newcommand{\bs}[1]{\boldsymbol{#1}}
\renewcommand{\bf}[1]{\mathbf{#1}}
\renewcommand{\rm}[1]{\mathrm{#1}}

\newcommand{\bbN}{\mathbb{N}}

\newcommand{\bbZ}{\mathbb{Z}}

\newcommand{\sfE}{\mathsf{E}}

\newcommand{\A}{\mathcal{A}}
\newcommand{\B}{\mathcal{B}}
\newcommand{\C}{\mathcal{C}}

\newcommand{\F}{\mathcal{F}}

\renewcommand{\L}{\mathcal{L}}
\renewcommand{\P}{\mathcal{P}}
\newcommand{\Q}{\mathcal{Q}}
\newcommand{\R}{\mathcal{R}}

\newcommand{\U}{\mathcal{U}}
\newcommand{\V}{\mathcal{V}}
\newcommand{\W}{\mathcal{W}}
\newcommand{\Y}{\mathcal{Y}}

\newcommand{\frX}{\mathfrak{X}}

\renewcommand{\S}{\Sigma}

\renewcommand{\a}{\alpha}
\renewcommand{\b}{\beta}
\newcommand{\eps}{\varepsilon}

\renewcommand{\k}{\kappa}
\renewcommand{\l}{\lambda}

\newcommand{\s}{\sigma}

\newcommand{\spt}{\mathrm{spt}}

\newcommand{\ol}[1]{\overline{#1}}

\newcommand{\into}{\hookrightarrow}
\newcommand{\onto}{\twoheadrightarrow}
\newcommand{\fin}{\nolinebreak\hspace{\stretch{1}}$\lhd$}

\newcommand{\llangle}{\langle\!\langle}
\newcommand{\rrangle}{\rangle\!\rangle}
\renewcommand{\t}[1]{\tilde{#1}}

\begin{document}

\title{\textbf{Amenable groups with very poor compression into Lebesgue spaces}}
\author{Tim Austin}
\date{}

\maketitle

\begin{abstract}
We give a construction of finitely-generated amenable groups that
do not admit any coarse $1$-Lipschitz embedding with positive
compression exponent into $L_p$ for any $1 \leq p < \infty$, including some that are four-step solvable, answering positively a question of Arzhantseva, Guba and Sapir.
\end{abstract}

\parskip 0pt

\tableofcontents

\parskip 7pt

\section{Introduction}

Given a metric space $(X,\rho)$ and a Banach space $\frX$ with norm
$\|\cdot\|$, the \textbf{compression exponent}
$\a^\ast_\frX(X,\rho)$ of $X$ into $\frX$ is the supremum of those
$\a \geq 0$ for which there exists an injection $f:X\into \frX$ such
that
\[\rho(x,y)^\a \lesssim \|f(x) - f(y)\| \leq \rho(x,y)\quad\quad\forall x,y \in X\]
(where, as usual, $\lesssim$ and $\gtrsim$ denote inequalities that
hold up to an arbitrary positive multiplicative constant that is
independent of the arguments of the functions in question).

This can be viewed as a quantitative measure of how well, if at all,
the space $(X,\rho)$ can be coarsely Lipschitzly embedded into
$\frX$. It was introduced by Guentner and Kaminker
in~\cite{GueKam04}, with particular emphasis on the case when
$(X,\rho)$ is a finitely-generated group equipped with a
left-invariant word metric associated to a finite symmetric
generating set, and $\frX$ is a Hilbert space. In this case the
value $\a_\frX^\ast(G)$ is easily seen to be independent of the
particular choice of generating set, since different generating sets
lead to quasi-isometric word metrics.  For a group $G$ we can define
analogously to $\a^\ast_\frX(G)$ the \textbf{equivariant compression
exponent} $\a^\#_\frX(G)$ as the supremum of those $\a \geq 0$ for
which there exists an \emph{equivariant} injection $f:G\into \frX$
(that is, an injection given by the orbit of a point under some
action $G\curvearrowright \frX$ by affine isometries) such that
\[\rho(g,h)^\a \lesssim \|f(g) - f(h)\| \leq \rho(g,h)\quad\quad\forall g,h \in G.\]
We write $\a^\ast_p$ and $\a^\#_p$ in place of $\a^\ast_{L_p}$ and
$\a^\#_{L_p}$. Since Guentner and Kaminker's work a number of
methods have been brought to bear on the estimation of these
exponents for different groups; see, for
instance,~\cite{ArzDruSap09,ArzGubSap06,deCStaVal08,deCTesVal07,Tes07,Tes08,Tes08.2,Tes09,AusNaoPer--ZwrZ,NaoPer08,NaoPer09}.
A more detailed discussion of these developments and more complete
references can be found in the introduction to~\cite{NaoPer09}.  It
is also worth noting that if $G$ is amenable then necessarily
$\a_2^\#(G) = \a_2^\ast(G)$, as was shown by Aharoni, Maurey and
Mityagin in~\cite{AhaMauMit85} for Abelian $G$ and then by Gromov
(see~\cite{deCTesVal07}) in general.

Most known results concern the case of amenable groups. It was shown
by Guentner and Kaminker in~\cite{GueKam04} that if $\a^\#_2(G) >
\frac{1}{2}$ then this actually implies that $G$ is amenable, and
another proof of this fact has now been given in~\cite{NaoPer08}
along with a generalization to other Lebesgue target spaces
(although it is known that there are both amenable and non-amenable
groups that have Euclidean equivariant compression exponent exactly
equal to $\frac{1}{2}$).  In the reverse direction, Arzhantseva,
Guba and Sapir asked as Question 1.12 in~\cite{ArzGubSap06} whether
there are amenable groups that have compression exponent $0$ (or,
more modestly, strictly less than $\frac{1}{2}$) for embeddings into
Hilbert space, and following progress in other directions versions
of this question have since been re-posed as Question 5.23 of
Arzhantseva, Drutu and Sapir~\cite{ArzDruSap09} and Question 1.5 of
Tessera~\cite{Tes09}.

This question is interesting partly because previous methods for
bounding Hilbert space compression exponents from above seem unable
in principle to break the $\frac{1}{2}$-barrier. In particular, a
very general upper bound in terms of random walk escape speed, first
introduced in~\cite{AusNaoPer--ZwrZ} and then considerably
generalized in~\cite{NaoPer08,NaoPer09}, cannot push below this
value. One version of that result (although not quite the most
general) asserts that if $\frX$ is a Banach space having modulus of
smoothness of power type $p$, then $\a_\frX^\#(G) \leq
\frac{1}{p\beta^\ast(G)}$, where $\beta^\ast(G)$ is the supremum of
those $\b \geq 0$ for which
\[\sfE(\rho(e_G,X_t))\gtrsim t^\b\quad\quad\forall t \in \bbN,\]
where $(X_t)_{t\geq 0}$ is the symmetric random walk on $G$ starting
at the identity $e_G$ corresponding to some finite generating set.
Since it is known that $\frac{1}{2} \leq \b^\ast \leq 1$ always,
this bound cannot give a value below $\frac{1}{2}$.

Moreover, in several concrete cases (in particular among certain
iterated wreath products of cyclic groups) the random walk bound
turns out to be the correct value of the compression exponent, even
when this value lies in $(\frac{1}{2},1)$.  These observations led
Naor and Peres to ask in~\cite{NaoPer09} (Question 10.3) whether it
is always the case that an amenable group $G$ has compression
exponent into $L_p$ given by exactly
$\min\{\frac{1}{p\beta^\ast},1\}$ (which would, in particular,
answer negatively the question of Arzhantseva, Guba and Sapir).

In this work we decide these questions by proving the following.

\begin{thm}\label{thm:main}
There is a finitely-generated amenable group that does not admit any
embedding into $L_p$ with a positive compression exponent for any $p
\in [1,\infty)$.
\end{thm}

Intuitively, this means we have found a finitely-generated amenable
group with much worse embedding properties into $L_p$ than have been
witnessed for such groups heretofore. Various notions relating to
the uniform embeddability of finitely-generated groups into `nice'
Banach spaces (predominantly Hilbert space) have been introduced in
geometric group theory (see, in particular, Gromov's discussion of
a-T-menability and some relatives in Section 7.E
of~\cite{Gro93}), and some equivalences found among them; on the
other hand, some of these properties have been shown to have
striking consequences elsewhere in gometric group theory and
topology (perhaps most notably Yu's deduction of the coarse
Baum-Connes and Novikov Conjectures therefrom in~\cite{Yu00}).
Gromov asked whether amenable groups are always a-T-menable, and
this was then proved by Bekka, Cherix and Valette
in~\cite{BekCheVal93}.  Thus from the relatively soft viewpoint of
coarse geometry amenable groups are known to possess such desirable
embeddings; the present paper indicates that this is \emph{not} a
guarantee of any very strong quantitative versions of the same
conclusions.

Before launching into technical details, we offer a sketch of our
approach to Theorem~\ref{thm:main}.

Underpinning the proof is a simple observation about how a sequence
of finite metric spaces with growing $L_p$-distortion can serve as
an obstruction to the good-compression $L_p$-embedding of an
infinite metric space. If $(X,\rho)$ is an infinite metric space,
and inside it we can find a family of bi-Lipschitzly embedded finite
metric spaces $(Y_n,\s_n)$, say with embeddings $\varphi_n:Y_n \into
X$, then these give rise to a bound on $\a_p^\ast(X,\rho)$ in terms
of
\begin{itemize}
\item the distortions $c_p(Y_n,\s_n)$, and
\item the ratios according to which distances under $\s_n$ are
approximately expanded by $\varphi_n$.
\end{itemize}
This interplay between the distortions of the $Y_n$ and approximate
expansion ratios of the $\varphi_n$ is crucial: the resulting
obstruction to $L_p$-embedding of the whole of $(X,\rho)$ will be
more severe according as the expansion ratios of these $\varphi_n$
grow more slowly.  A quantitative version of this observation that
is suitably tailored for our proof of Theorem~\ref{thm:main} appears
as Lemma~\ref{lem:obs} below. With this principle in hand, we will
find the specific finite metric spaces that will serve as our
obstructions in the form of certain quotients of high-dimensional
Hamming cubes. Regarding Hamming cubes as vector spaces over
$\bbZ_2$, a result of Khot and Naor~\cite{KhoNao06} provides us with
certain quotients of these by $\bbZ_2$-subspaces whose Euclidean
distortion is linearly large in their diameter.

Most of our work will then go into constructing a group that contains suitable bi-Lipschitz copies of these cube-quotients.  The basic strategy of using a sequence of poor-distortion finite subsets to bound above the compression exponent of a group is not new; for example, Arzhantzeva, Dru\c{t}u and Sapir use it in~\cite{ArzDruSap09} for their construction of finitely-generated groups with arbitrary compression exponents into uniformly convex Banach spaces, except that their finitary obstructions are sequences of expanders, and they use free products (hence giving non-amenable examples) to construct finitely-generated groups containing them.

Our construction will give a two-fold Abelian extension of any suitable
base group $G$ (which could be, for example, the classical
lamplighter $\bbZ_2\wr \bbZ$, so that among our examples of groups
with very poor compression we find certain four-step solvable
groups). The first of these extensions of $G$ is a wreath product,
but the second is something rather more complicated.

More precisely, we begin with any finitely-generated amenable $G$
that has exponential growth, form the wreath product over it $H :=
\bbZ_6\wr G$ and equip it with its natural lifted generating set and
word metric (we will comment on the modulus $6$ shortly). Consider $H$ as acting by translation on the
$\bbZ_2$-vector space $\bbZ_2^{\oplus H}$. These are the ingredients
needed to form another wreath product $\bbZ_2\wr H = \bbZ_2^{\oplus
H}\rtimes H$, but we will instead examine a slightly more
complicated relative of this construction. (Indeed, while the random
walk method gives an upper bound $\a_2^\ast(\bbZ_2\wr H) \leq
\frac{1}{2}$, and a recent analysis by Li~\cite{Li09} building on a
construction from de Cornulier, Stalder and
Valette~\cite{deCStaVal09} has shown that it is at least
$\frac{1}{6}$, its exact value is currently unknown. Arzhantseva,
Guba and Sapir originally proposed the related iterated wreath
product $\bbZ\wr (\bbZ\wr \bbZ)$ as a candidate for having zero
Euclidean compression exponent, but Li's work now also gives a
positive lower bound for this exponent.) We will first identify a
translation invariant $\bbZ_2$-subspace $V \leq \bbZ_2^{\oplus H}$,
and will then form instead the semidirect product $(\bbZ_2^{\oplus
H}/V)\rtimes H$, equipped with the natural generating set that
consists of lifts of the generators of $H$ and the single new
generator $\delta_{e_H} + V \leq \bbZ_2^{\oplus H}/V$.

The r\^ole of this $\bbZ^2$-subspace $V$ is to insert copies of a sequence of Khot and Naor's poor-distortion cube-quotients, say $\bbZ_2^{I_n}/C_n$ for some increasingly large index-sets $I_1$, $I_2$, \ldots, into the `zero-section' subgroup
\[\bbZ_2^{\oplus H}/V \subset (\bbZ_2^{\oplus H}/V)\rtimes H.\]
More precisely, we will prove that for some sequence of these cube-quotients $\bbZ_2^{I_n}/C_n$, one can effectively choose embeddings $\varphi^\circ_n:\bbZ_2^{I_n} \into \bbZ_2^{\oplus H}$ and a translation-invariant subgroup $V\leq \bbZ_2^{\oplus H}$ so that $\varphi^\circ_n(C_n)\subseteq V$, and the resulting quotient maps
\[\varphi_n:\bbZ_2^{I_n}/C_n \into \bbZ_2^{\oplus H}/V \subset (\bbZ_2^{\oplus H}/V)\rtimes H\]
are bi-Lipschitz embeddings with control on their distortions and expansion ratios as required for them to constrain the compression exponents of $L_p$-embeddings of the whole group $(\bbZ_2^{\oplus H}/V)\rtimes H$.

The principal difficulties in this construction arise from the need to insert a whole infinite sequence of these cube-quotients $\bbZ_2^{I_n}/C_n$ into $\bbZ_2^{\oplus H}/V$ using a single choice of subgroup $V$, with a uniform bound on their distortions and slow growth of their expansion ratios.  It is to address these difficulties that we specifically take $H = \bbZ_6\wr G$ for some $G$ of exponential growth.  The $\bbZ_2$-subspace $V$ will be assembled as an infinite sum $V_1 + V_2 + \ldots$ of $\bbZ_2$-subspaces, each of which creates the corresponding embedding of $\bbZ_2^{I_n}/C_n$, and we will use the particular algebraic structure of this wreath product $H$ to choose these $V_n$ so that they do not `interact' with each other too much: more specifically, so that quotienting by $V_{n'}$ for $n' > n$ does not disrupt the bi-Lipschitz embedding of $\bbZ_2^{I_n}/C_n$ that we introduce upon quotienting by $V_n$.  One feature of this algebraic argument is that it leads to the careful choice of $6$ (rather than, say, $2$ or $4$) as the modulus in the wreath product $H = \bbZ_6\wr G$.  This is used via the fact that there exists a non-constant function $\A:\bbZ_6\to\bbZ_2$ with the property that one cannot form the constant function $1 \in \bbZ_2$ as a $\bbZ_2$-linear comination of translates of $\A$ by elements of the domain $\bbZ_6$.  For this we take $\A$ to be the $\bbZ_2$-valued indicator function of $\{0,1,3,4\}$.  Since at a different step in the construction we will want a modulus that is even, a quick check shows that $6$ is the smallest possibility. By allowing $\A$ to act on different coordinates in $\bbZ_6^{\oplus G}$ we can create a large family of functions $\bbZ_6^{\oplus G}\to \bbZ_2$ (that is, elements of $\bbZ_2^{\oplus \bbZ_6^{\oplus G}}$) all of whose translates are linearly independent, and this will then be used in the proof that our choices of translation invariant subspaces $V_1,V_2,\ldots \leq \bbZ_2^{\bbZ_6\wr G}$ are `well-separated' (see Lemma~\ref{lem:lin-indep-+} and its use in the proof of Proposition~\ref{prop:TS-cubes}).

\textbf{Acknowledgements}\quad I am grateful to Assaf Naor and Yuval
Peres for introducing me to the topics touched by this paper and for
many subsequent discussions, to Sean Li and Assaf Naor for each
suggesting a number of improvements to this paper, and to Microsoft
Research Redmond for a period of their hospitality during which most
of this work was completed. \fin

\section{Background and generalities}

\subsection{Metric spaces}

The following nomenclature distinguishing certain kinds of metric
space will prove useful.

\begin{dfn}
A metric space $(X,\rho)$ is
\begin{itemize}
\item \textbf{infinite} if this is so of $X$ as a set;
\item \textbf{locally finite} if every open ball $B_X(x,r)$ in $X$ is a
finite set;
\item \textbf{$d_{\rm{min}}$-discrete} for some $d_{\rm{min}} > 0$ if
$\rho(x,y)\geq d_{\rm{min}}$ whenever $x,y \in X$ are distinct.
\end{itemize}
\end{dfn}

Given two metric spaces $(X,\rho)$ and $(Z,\theta)$ and a map
$f:X\to Z$, the \textbf{distortion} of $f$ is a measure of the
extent to which $f$ fails to be a homothety:
\[\rm{distortion}(f) := \sup_{u,v\in X,\,u\neq v}\frac{\theta(f(u),f(v))}{\rho(u,v)}\cdot\sup_{u,v\in X,\,u\neq v}\frac{\rho(u,v)}{\theta(f(u),f(v))},\]
where this is interpreted as $+\infty$ unless $f$ is at least a
bi-Lipschitz embedding.  Relatedly, the map $f$ is \textbf{$M$-bi-Lipschitz} for some $M \geq 1$ if it is bi-Lipschitz with distortion at most $M$.

The \textbf{distortion of $X$ into $Z$} is
now obtained by infimizing over $f$: it is conventionally denoted by
\[c_{(Z,\theta)}(X,\rho) := \inf_{f:X\to Z}\rm{distortion}(f).\]

The case in which $(Z,\theta)$ is a Lebesgue space $L_p$ with its
norm metric for some $1 \leq p < \infty$ is particularly
well-studied, and in this case we abbreviate $c_{L_p}$ to $c_p$.

In addition to the distortion of a map $f$, we will sometimes need
to keep track of its \textbf{expansion ratio}: in the above notation
this is simply the ratio
\[\inf_{u,v\in X,\,u\neq v}\frac{\theta(f(u),f(v))}{\rho(u,v)}.\]
If the expansion ratio is $r$, then the above definitions combined
tell us that
\[r\cdot \rho(u,v) \leq \theta(f(u),f(v))\leq \rm{distortion}(f)\cdot r\cdot \rho(u,v)\quad\quad\forall u,v \in X,\]
provided all the quantities appearing here are finite.

An important class of metrics that will appear repeatedly later is that of the Hamming metrics:

\begin{dfn}[Hamming metric]
If $(X,\rho)$ is a metric space and $n\geq 1$ then the \textbf{Hamming metric} on $X^n$ associated to $\rho$ is defined by
\[d_{\rm{Ham}}\big((x_i)_{i\leq n},(y_i)_{i\leq n}\big) := \sum_{i=1}^n\rho(x_i,y_i).\]
\end{dfn}

We will also make use of the construction of quotient metrics.  The
following definition is discussed, for example, in Section 3 of Khot
and Naor~\cite{KhoNao06}.

\begin{dfn}[Quotient metrics]
Suppose that $(X,\rho)$ is a metric space and that $\P$ is a
partition of $X$ with the property that
\begin{quote}
For any $P,Q \in \P$ and any $x \in P$ there is some $y \in Q$ such
that $\rho(x,y) = \rm{dist}_{\rho}(P,Q)$ (that is, `the minimal
distance between $P$ and $Q$ can be realized from any starting point
on one side').
\end{quote}

Then the \textbf{quotient metric} $\rho_{/\P}$ on $\P$ is defined by
\[\rho_{/\P}(P,Q) := \rm{dist}_{\rho}(P,Q).\]
The routine verification that this is indeed a metric on $\P$ under
the above assumption can be found in Section 3 of~\cite{KhoNao06}.

In particular, this definition always applies to give a metric on
the homogeneous space $G/H$ when $(G,\rho)$ is a group carrying a
left-invariant metric $\rho$ and $H \leq G$.  In this case we denote
the quotient metric by $\rho_{/H}$.
\end{dfn}

Before leaving this subsection, let us recall some useful analyst's notation.  Given positive quantities $A$ and $C$ and any other parameter $B$, we write $A \lesssim_B C$ or $A = \rm{O}_B(C)$ if quantity $A$ is bounded above by quantity $C$ up to a positive multiplicative constant depending only on $B$ (and simply $A \lesssim C$ or $A = \rm{O}(C)$ if the constant is truly universal); similarly $A = \Omega_B(C)$ if $A$ is bounded from below by such a positive multiple of $C$; and $A = \Theta(C)$ if $A$ is bounded from above and below by the multiples of $C$ by two universal positive bounds (that is, if $A = \rm{O}(C)$ and $A = \Omega(C)$).  Thus, for instance, an embedding has `expansion ratio $\Theta(1)$' if its expansion can be shown to lie in $[r,R]$ for some $R \geq r > 0$, independently of any other parameters that went into its construction.

\subsection{Groups}

The construction of this paper makes a double use of a semidirect
product.  Given a group $L$ and another $H$ that acts on it by
automorphisms, $H\curvearrowright L$, we may form the
\textbf{semidirect product} $H \ltimes L$ or $L\rtimes H$ as the set
$L\times H$ equipped with the group operation
\[(\ell_1,h_1)\cdot(\ell_2,h_2) := (\ell_1^{h_2}\cdot\ell_2,h_1h_2).\]

Importantly, the group $H\ltimes L$ is finitely-generated if
\begin{itemize}
\item $H$ is finitely-generated and
\item $L$ admits a finite subset $S \subseteq L$ such that $L =
\langle \{\a^h(s):\ h \in H,\, s \in S\}\rangle$.
\end{itemize}
This can be so even if $L$ itself is \emph{not} finitely-generated.
An important example of this is the wreath product $K \wr H :=
K^{\oplus H}\rtimes H$, formed from the action $H\curvearrowright
K^{\oplus H}$ (the group of $H$-indexed families of members of $K$
with cofinitely many entries equal to $e_K$) by coordinate
right-translation so that
\begin{eqnarray}\label{eqn:semidirect}
\big((k_{1,h})_{h\in H},h_1\big)\cdot \big((k_{2,h})_{h\in
H},h_2\big) &=& \big((k_{1,hh_2^{-1}}k_{2,h})_{h\in H},h_1h_2\big).
\end{eqnarray}
This is always finitely-generated if this is so of $H$ and $K$, even
though if $H$ is infinite then $L = K^{\oplus H}$ is not
finitely-generated.

The semidirect product $L\rtimes H$ always contains a canonical copy
of $L$ in the form of the subgroup $\{(\ell,e_H):\ \ell \in L\}$:
following Naor and Peres~\cite{NaoPer09} we refer to this as the
\textbf{zero section} of $L\rtimes H$ and sometimes denote it by
$(L\rtimes H)_0$.

In this paper we will use a construction similar to that of the
wreath product, except that the place of $K^{\oplus H}$ will
sometimes be taken by one of its quotients. In the following the
group $K$ will always be one of the cyclic groups $\bbZ_2$ or
$\bbZ_6$.  We will later form first an extension of a suitable base
group by a power of $\bbZ_6$, and then a further extension by a
power of $\bbZ_2$, and our choice of notation is geared to
distinguish between these two extensions as clearly as possible. In
particular:
\begin{itemize}
\item Given a base group $G$, we will usually denote elements of $\bbZ_6^{\oplus
G}$ by lowercase bold letters such as $\bf{w} = (w_g)_{g\in G}$, and
refer to them as \textbf{vectors}.  In this setting we write
$\bf{e}_g$ for $g \in G$ for the vector with $g^{\rm{th}}$ entry
equal to $1 \in \bbZ_6$ and all other entries equal to $0$; we refer
to these $\bf{e}_g$ collectively as the \textbf{standard generators}
of $\bbZ_6^{\oplus G}$.
\item On the other hand, given another base group $H$ (which will later be equal to $\bbZ_6\wr
G$), we will denote elements of $\bbZ_2^{\oplus H}$ by uppercase
calligraphic letters, and will refer to them as \textbf{functions}
$\W:H \to \bbZ_2$; correspondingly an expression such as $\W + \V$ refers to a sum of functions defined pointwise. In this setting we denote by $\delta_h$ the
function that takes the value $1\in \bbZ_2$ at $h$ and $0$
elsewhere. We will sometimes need to work instead with the set $\{h
\in H:\ \W(h) = 1\}$, which we refer to as the \textbf{support} of
$\W$ and denote by $\spt\,\W$.
\end{itemize}
The only exceptions to these rules are that $\bs{0}$ will be used
for the zero element in either case, since this should cause no
confusion, and that for any finite index set $T$ we write
$\bs{1}_{T}$ for the function $T\to C$ that identically takes the
value $1 \in C$ for either $C = \bbZ_2$ or $C = \bbZ_6$.

If $V \leq \bbZ_2^{\oplus H}$ is a subspace that is invariant under
the coordinate right-translation action of $H$, then this action of
$H$ quotients to a well-defined action $H \curvearrowright
\bbZ_2^{\oplus H}/V$.  Given also a chosen finite symmetric
generating set $S$ for $H$, we will always endow the semidirect
product $(\bbZ_2^{\oplus H}/V)\rtimes H$ with the symmetric
generating set
\[\{(\delta_{e_H} + V,e_H)\}\cup \{(\bs{0},s):\ s\in S\}.\]
To this enlarged generating set we always associate the
left-invariant word metric as usual.  In general, if the generating
set $S$ of $H$ is understood, the above generating set will also be
understood for $(\bbZ_2^{\oplus H}/ V)\rtimes H$, and we will
(slightly abusively) denote the resulting metric by
$\rho_{(\bbZ_2^{\oplus H}/V)\rtimes H}$, and will denote the
restriction of this new metric to the zero section by
$\rho_{((\bbZ_2^{\oplus H}/V)\rtimes H)_0}$.

We will be working largely with the restriction of this word metric
to the zero section of $(\bbZ_2^{\oplus H}/ V)\rtimes H$, and to
this end it will be helpful to have a simpler equivalent metric on
the zero section to work with.  Such a description can quite easily
be given in terms of traveling salesman tours as a simple extension
of the usual heuristic interpretation of the wreath product word
metric in terms of lamplighter walks. The definition below has been
adapted from Section 2 of Naor and Peres~\cite{NaoPer09}, where its
connection with the problem of estimating the length of traveling
salesman tours among points in a metric space (and, in particular,
with the work of Jones in~\cite{Jon90}) is explained in more detail.

\begin{dfn}[Pinned traveling salesman metrics]\label{dfn:TS}
Given a metric space $(X,\rho)$ and a distinguished point $x^\circ
\in X$, the associated \textbf{traveling salesman metric pinned at
$x_0$} is the metric $\rm{TS}_{\rho,x^\circ}$ on the collection
$\bbZ_2^{\oplus X}$ of finitely-supported maps $X\to \bbZ_2$
(equivalently, finite subsets of $X$) defined by
\begin{multline*}
\rm{TS}_{\rho,x^\circ}(\A,\B)\\ \quad\quad:=
\left\{\begin{array}{ll}0&\quad\hbox{if }\A = \B\\
\min\left\{\begin{array}{l}\sum_{i=1}^\ell\rho(x_i,x_{i+1}):\\
\quad\hbox{the cycle }(x_1 = x^\circ,x_2,\ldots,x_\ell,x_{\ell + 1}
= x^\circ)\\ \quad\hbox{ covers }\spt(\A + \B)\end{array}\right\} +
1&\quad\hbox{else.}\end{array}\right.
\end{multline*}
This is clearly a metric in view of the inclusion $\spt(\A +
\C)\subseteq \spt(\A + \B) \cup \spt(\B + \C)$ and the ability to
concatenate covering cycles.

Abbreviating, it will be understood that the notation
$\rm{TS}_{\rm{Ham}}$ on $\bbZ_2^{\oplus \bbZ_6^{\oplus S}}$ refers
to the Hamming metric $d_{\rm{Ham}}$ on the underlying group
$\bbZ_6^{\oplus S}$ (where $\bbZ_6$ is given the word metric
corresponding to the generators $\{1,-1\}$) and the distinguished
point $x^\circ := \bs{0} \in \bbZ_6^{\oplus S}$.

If $T$ is a finite index set and $F \leq \bbZ_2^{\bbZ_6^T}$ is a
$\bbZ_2$-subspace invariant under the translation action of
$\bbZ_6^T$, then we denote by $\rm{TS}_{\rm{Ham}/F}$ the metric on
$\bbZ_2^{\bbZ_6^T}/F$ that results from quotienting
$\rm{TS}_{\rm{Ham}}$.
\end{dfn}

Let us note at once that following trivial consequence of this
definition:

\begin{lem}\label{lem:monotone}
If $\V,\V',\W,\W':X \to \bbZ_2$ are finitely supported and $x^\circ
\in X$ then
\[\spt\,(\V' + \W') \supseteq \spt\,(\V + \W)\quad\quad\Rightarrow\quad\quad \rm{TS}_{\rho,x^\circ}(\V',\W')\geq \rm{TS}_{\rho,x^\circ}(\V,\W).\]
In particular, $\rm{TS}_{\rho,x^\circ}$ is addition-invariant on
$\bbZ_2^{\oplus X}$. \qed
\end{lem}

A more delicate calculation that we will need later is the
following.  This closely resembles the derivation of equation (13)
in Naor and Peres~\cite{NaoPer09}, and we offer only a rather terse
account here.

\begin{lem}\label{lem:approx-by-TS}
If $(\bbZ_2^{\oplus H}/V)\rtimes H$ is a semidirect product as
above, $H$ has generating set $S$ and corresponding word metric
$\rho$ and we lift this to a generating set of $(\bbZ_2^{\oplus
H}/V)\rtimes H$ as above with associated metric
$\rho_{(\bbZ_2^{\oplus H}/V)\rtimes H}$, then the inclusion
$\bbZ_2^{\oplus H}/V \equiv ((\bbZ_2^{\oplus H}/V)\rtimes H)_0
\subset (\bbZ_2^{\oplus H}/V)\rtimes H$ is a $2$-bi-Lipschitz
embedding of $\rm{TS}_{\rho,e_H\,/V}$ into $\rho_{(\bbZ_2^{\oplus
H}/V)\rtimes H}$ with expansion ratio lying in $\Theta(1)$.
\end{lem}

To find a shortest-length word that evaluates to a given group element $(\A + V,e_H) \in (\bbZ_2^{\oplus H}/V)\rtimes H$ amounts to finding a shortest walk from $(V,e_H)$ to $(\A + V,e_H)$ that steps only along edges of the Cayley graph defined by the generated set of $(\bbZ_2^{\oplus H}/V)\rtimes H$.  Intuitively, in order to do this we must take a walk in the Cayley graph of $H$ that starts and ends at $e_H$, and along the way pass through finitely many points of $H$ at which we modify the value taken by the identically-zero function $\bs{0}:H\to \bbZ_2$ to obtain a function that lies in $\A + V$ (this much corresponds roughly to an evaluation of a distance in $\rm{TS}_{\rho,e_H}$), and finally infimize over the choice of that element of $\A + V$ (which corresponds to working in the quotient metric $\rm{TS}_{\rho,e_H\,/V}$).  The following proof makes this intuition precise.

\textbf{Proof}\quad We will show that for any finitely-supported
$\A,\B:H \to \bbZ_2$ we have
\begin{multline*}
\rm{TS}_{\rho,e_H\,/V}(\A + V,\B + V) \leq \rho_{(\bbZ_2^{\oplus
H}/V)\rtimes H}((\A + V,e_H),(\B + V,e_H))\\ \leq 2\cdot
\rm{TS}_{\rho,e_H\,/V}(\A + V,\B + V).
\end{multline*}

To see this, recall that by definition $\rho_{(\bbZ_2^{\oplus
H}/V)\rtimes H}((\A + V,e_H),(\B + V,e_H))$ is the shortest length
of any word in the alphabet given by the generating set
\[\{(\delta_{e_H} + V,e_H)\} \cup \{(V,s):\ s\in S\}\]
whose evaluation in $(\bbZ_2^{\oplus H}/V)\rtimes H$ is equal to
$(\A + V,e_H)^{-1}\cdot (\B + V,e_H) = (\B - \A + V,e_H)$.  Set
$\ell := \rho_{(\bbZ_2^{\oplus H}/V)\rtimes H}((\A + V,e_H),(\B +
V,e_H))$, let $g_\ell g_{\ell - 1}\cdots g_1$ be such a word of
minimal length in this alphabet and define $G_i := g_ig_{i-1}\cdots
g_1 =: (\C_i + V,h_i)$ for $i=1,2,\ldots,\ell$ and also $G_0 :=
(V,e_H) =: (\C_0 + V,h_0)$.

Each $G_i$ differs from $G_{i-1}$ by left-multiplication by either
$(\delta_{e_H} + V,e_H)$ or by $(V,s)$ for some $s\in S$.  In the
first case the multiplication rule~(\ref{eqn:semidirect}) tells us
that $h_i = h_{i-1}$ and $\C_i = \C_{i-1} + \delta_{h_{i-1}} \mod
V$, and in the second it gives $\C_i = \C_{i-1}$ but $h_i =
g_ih_{i-1}$.  Letting $1 \leq i_1 < i_2 < \ldots < i_k \leq \ell$ be
the subsequence of those $i$ at which we are in the first of these
situations, we see that overall
\[\C_\ell = \delta_{h_{i_1-1}} + \delta_{h_{i_2-1}} + \cdots + \delta_{h_{i_k-1}} = \A - \B\mod V.\]
Therefore, omitting the steps $i_1$, $i_2$, \ldots, $i_k$, the
sequence $h_1$, $h_2$, \ldots, $h_{i_1 - 1}$, $h_{i_1 + 1}$, \ldots,
$h_\ell$ executes a walk in $H$ that starts from $e_H$, visits every
point of $\spt\,\C_\ell$ (where $\C_\ell \in \A - \B + V$) at least
once and then returns to $e_H$; and at those omitted steps the
multiplication by $g_{i_j}$ does not change the position of $h_{i_j
- 1}$ but instead alters the function $\C_{i_j-1}$ to $\C_{i_j}$
modulo $V$.

On the one hand, it follows immediately from Definition~\ref{dfn:TS}
that the length of this walk is at least
$\rm{TS}_{\rho,e_H}(\bs{0},\C_\ell)$; and on the other, for any such
$\C_\ell$, since any two points of $\spt\,\C_\ell$ are separated by
a distance of at least $1$ in $H$ under $\rho$, it follows that
there does exist such a sequence of walk-steps interspersed with
additions of $\delta_{h_{i_k}} + V$ of length at most
\[\rm{TS}_{\rho,e_H}(\bs{0},\C_\ell) + |\C_\ell| \leq 2\cdot \rm{TS}_{\rho,e_H}(\bs{0},\C_\ell).\]

Therefore, first minimizing the above lower bound over the possible
choices of $\C_\ell \in \A - \B + V$ we obtain
\begin{multline*}
\rm{TS}_{\rho,e_H/V}(\A + V,\B + V) = \min_{\C \in \A - \B +
V}\rm{TS}_{\rho,e_H}(\bs{0},\C)\\ \leq \rho_{(\bbZ_2^{\oplus
H}/V)\rtimes H}((\A + V,e_H),(\B + V,e_H)),
\end{multline*}
and secondly choosing $\C = \C_\ell$ that attains this minimum and
then selecting a suitable word $g_\ell g_{\ell -1}\cdots g_1$ as
described above proves the corresponding upper bound. \qed

\section{Using distortion to bound compression}

Before turning to the construction of our group, we examine the
relation between the notions of compression for infinite metric
spaces and of distortion for finite metric spaces.  In particular,
we will see how to use large high-distortion finite subsets of an
infinite metric space as obstructions to good-compression embeddings
of the whole space, in the sense made precise by the following.

\begin{lem}\label{lem:obs}
Suppose that $\frX$ is a normed vector space, that $(X,\rho)$ is an
infinite, locally finite, $1$-discrete metric space, and suppose
further that we can find a sequence of finite $1$-discrete metric
spaces $(Y_n,\s_n)$ and embeddings $\varphi_n:Y_n \into X$ such that
\begin{itemize}
\item the $Y_n$ are increasing in diameter:
$\rm{diam}(Y_n,\s_n)\to\infty$;
\item the $Y_n$ are embedded in $X$ with uniformly-bounded distortion: there
are some fixed $L \geq 1$ and some sequence of positive reals
$(r_n)_{n\geq 1}$ such that
\[\frac{1}{L}r_n\s_n(u,v) \leq \rho(\varphi_n(u),\varphi_n(v)) \leq Lr_n\s_n(u,v)\quad\quad\forall u,v \in Y_n,\,n\geq 1\]
(so the distortion is at most $L^2$);
\item the $Y_n$ do not expand too fast inside $X$ relative to their
size: we have $1 \leq r_n \lesssim \rm{diam}(Y_n,\s_n)^\eps$ for
every $\eps > 0$;
\item the $Y_n$ have bad distortion into $\frX$: for some $\eta >
0$ we have $c_\frX(Y_n,\s_n) \gtrsim \rm{diam}(Y_n,\s_n)^\eta$ for
all $n \geq 1$.
\end{itemize}
Then $\a_\frX^\ast(X,\rho) \leq 1-\eta$.
\end{lem}

\textbf{Proof}\quad We may assume that $\a^\ast(X,\rho) > 0$, since
otherwise the result is trivial.  Given this, let $\a <
\a^\ast(X,d)$ and let $f:X\to\frX$ be a $1$-Lipschitz embedding into
$\frX$ that achieves compression $\a$:
\[\rho(x,y)^\a \lesssim \|f(x) - f(y)\| \leq \rho(x,y)\quad\quad \forall x,y \in X.\]

Combining this with the bi-Lipschitz condition
\[\frac{1}{L}r_n \s_n(u,v) \leq \rho(\varphi_n(u),\varphi_n(v)) \leq Lr_n\s_n(u,v)\]
by setting $x := \varphi_n(u)$ and $y := \varphi_n(v)$ in the
compression inequality, we deduce that $f\circ \varphi_n: Y_n\into
\frX$ satisfies
\[\frac{1}{L^\a}\big(r_n\s_n(u,v)\big)^\a \lesssim \|f\circ \varphi_n(u) - f\circ \varphi_n(v)\| \leq Lr_n\s_n(u,v)\quad\quad\forall u,v \in Y_n,\ n\geq 1.\]
It follows that
\begin{eqnarray*}
&&\rm{distortion}(f\circ \varphi_n)\\
&&\quad\stackrel{\rm{dfn}}{=} \max_{u,v \in Y_n,\,u \neq
v}\frac{\s_n(u,v)}{\|f\circ\varphi_n(u) -
f\circ\varphi_n(v)\|}\cdot\max_{u,v \in Y_n,\,u \neq
v}\frac{\|f\circ\varphi_n(u) - f\circ\varphi_n(v)\|}{\s_n(u,v)}\\
&&\quad\lesssim (\max_{u,v \in Y_n,\,u \neq v} L^\a
r_n^{-\a}\s_n(u,v)^{1
- \a})\cdot Lr_n\\
&&\quad= L^{1 + \a}r_n^{1-\a}\rm{diam}(Y_n,\s_n)^{1 - \a}.
\end{eqnarray*}

Since we know that this must be at least $c_\frX(Y_n,\s_n) \gtrsim
\rm{diam}(Y_n,\s_n)^\eta$, it follows that
\[\rm{diam}(Y_n,\s_n)^\eta \lesssim L^{1 + \a}r_n^{1-\a}\rm{diam}(Y_n,\s_n)^{1 - \a}.\]

Finally, from our assumption that $r_n \lesssim
\rm{diam}(Y_n,\s_n)^\eps$ for every $\eps > 0$ and since $L$ is
independent of $n$, it follows that
\[\rm{diam}(Y_n,\s_n)^\eta \lesssim \rm{diam}(Y_n,\s_n)^{(1 - \a)(1 + \eps)}\quad\quad\hbox{as}\ n\to\infty\]
for every $\eps > 0$, and hence that $1 - \a \geq \eta$, or $\a \leq
1 - \eta$. \qed

\textbf{Remark}\quad Note that the interplay between the expansion
ratio $r_n$ incurred by the embedding $\varphi_n:Y_n \into X$ and
the compression assumption on $f$ is crucial for this proof: if, in
a different situation, we knew that we could embed the same metric
spaces $Y_n$ into $X$ bi-Lipschitzly, but only at the expense of
drastically enlarging them, then the above argument would give no
bound on the compression. \fin

In conjunction with the general principle contained in the above
lemma, we will use the following result of Khot and
Naor~\cite{KhoNao06} to supply us with high-distortion finite metric
spaces to serve as obstructions.

\begin{thm}\label{thm:KhoNao}
There are $\bbZ_2$-subspaces $V_d \leq \bbZ_2^d$ for all $d \geq 1$
such that \[c_p(\bbZ_2^d/V_d,d_{\rm{Ham}/V_d}) \gtrsim_p d\] as $d
\to\infty$ for all $p \in [1,\infty)$. Moreover, the subspace $V_d$
may be taken of the form $C_d^\perp$ for some other subspace $C_d
\leq \bbZ_2^d$ such that $\bs{1} = (1,1,\ldots,1) \in C_d$, where
\[C_d^\perp := \{(w_i)_{i\leq d} \in \bbZ_2^d:\ w_1v_1 + w_2v_2 + \cdots + w_dv_d = 0 \in \bbZ_2\ \forall (v_i)_{i\leq d} \in C_d\}.\]
\qed
\end{thm}

\textbf{Proof}\quad The existence of such subspaces is covered
directly by Remark 3.1 in~\cite{KhoNao06}, and the fact that we may
take $\bs{1} \in C_d$ follows at once from the greedy-algorithm
proof given for their Corollary 3.5. \qed

With these ingredients to hand, we can now state a more precise
result from which Theorem~\ref{thm:main} will follow immediately.

\begin{thm}\label{thm:main2}
For any finitely-generated amenable group $G$ of exponential growth,
with finite symmetric generating set $S$ and associated
left-invariant word metric $\rho$, there are a fixed $M\geq 1$ and a
$\bbZ_2$-subspace $V \leq \bbZ_2^{\oplus (\bbZ_6\wr G)}$ such that
\begin{itemize}
\item $V$ is invariant under the action of $\bbZ_6\wr G$ by coordinate
right-translation, and
\item the natural choice of generating set for
$(\bbZ_2^{\oplus (\bbZ_6\wr G)}/V)\rtimes (\bbZ_6\wr G)$ admits a
sequence of injections
\[\varphi_n:\bbZ_2^{d_n}/V_{d_n} \into \bbZ_2^{\oplus (\bbZ_6\wr G)}/V\subset (\bbZ_2^{\oplus (\bbZ_6\wr G)}/V)\rtimes (\bbZ_6\wr G)\]
that are embeddings of $d_{\rm{Ham}/V_{d_n}}$ into
$\rho_{(\bbZ_2^{\oplus (\bbZ_6\wr G)}/V)\rtimes (\bbZ_6\wr G)}$ with uniformly-bounded distortions and expansion ratios $r_n \lesssim \log d_n$.
\end{itemize}
\end{thm}

\textbf{Proof of Theorem~\ref{thm:main} from
Theorem~\ref{thm:main2}}\quad This now follows by applying
Lemma~\ref{lem:obs} using the Euclidean distortion bound from Theorem~\ref{thm:KhoNao}, according to which we can take
$\eta = 1$ when the spaces $(Y_n,\s_n)$ are the cube-quotients
$(\bbZ_2^{d_n}/V_{d_n},d_{\rm{Ham}/V_{d_n}})$. \qed

The rest of this paper will be taken up by the proof of
Theorem~\ref{thm:main2}.  All the delicacy in this construction lies
in the selection of the subspace $V$. The formation of a tower of
two semidirect products (starting from a base group $G$ of
exponential growth, which could also be obtained as another
semidirect product) is rather important for the detailed arguments
to come, and leads to examples, at their simplest, that are $4$-step
solvable.  Note that they are, at least, elementary amenable. The
ways in which we use each of these semidirect products, and also the
reasons for choosing the particular cyclic groups $\bbZ_2$ and
$\bbZ_6$ as building blocks, will be examined as they arise during
the subsequent sections.

\section{Construction of the example}

\subsection{Approximately equidistant sets}

Our first step in the construction of the subspace $V$ promised by
Theorem~\ref{thm:main2} is the following.

\begin{lem}[Finding not-too-wide sets of many roughly equidistant
points]\label{lem:big-equidist} Let $G$ be a finitely-generated
group of exponential growth, $S$ a symmetric generating set and
$\rho$ the resulting word metric.  Then there are some $M \geq 1$
and $c > 1$, depending only on $G$ and $S$, such that for any $r_0 >
0$ there is some $r \geq r_0$ for which we can find a subset
$\{x_1,x_2,\ldots,x_d\}$ with $d := 10^r$ such that $r \leq
\rho(x_i,x_j) \leq 2Mr$ for all $1 \leq i < j \leq d$ and $r \leq
\rho(x_i,e_G) \leq 2Mr$ for all $1 \leq i \leq d$.
\end{lem}

\textbf{Proof}\quad Let $c > 1$ be such that
\[|B(e,r)| \geq c^r\quad\quad\hbox{for infinitely many}\ r \geq 1.\]
In the reverse direction, a na\"\i ve count tells us that
\[|B(e,r)| \leq |S|^r\quad\quad\forall r \geq 1.\]
Let $M \geq 1$ be a fixed integer such that $c^M > (11|S|)^2$. Then
for any $r_0 \geq 1$ we can find some $r \geq r_0 + 1$ and $0 \leq s
\leq M-1$ such that
\[|B(e,Mr)|\geq |B(e,M(r-1) + s)| \geq (c^M)^{r-1}\cdot c^s \geq (11|S|)^{2(r-1)} \geq 11^r\cdot |S|^r,\]
and the above inequalities together now imply that
\[|B(e,Mr)\setminus B(e,r)| = |B(e,M r)| - |B(e,r)| \geq (11^r
- 1)\cdot|S|^r \geq 10^r \cdot|S|^r.\]

Using again that $|B(x,r)| \leq |S|^r$ for every $(x,r)$, it follows
that $B(e,M r)\setminus B(e,r)$ cannot be covered by fewer than
$10^r$ balls of radius $r$. This implies instead that we can find in
$B(e,M r)\setminus B(e,r)$ a subset containing at least $d := 10^r$
points, any two of which are separated by a distance of at least
$r$. Enumerating this subset as $\{x_1,x_2,\ldots, x_d\}$ completes
the proof. \qed

\subsection{Some finite-dimensional quotient spaces}

By considering the subsets
\[\Big\{\Big(\sum_{i = 1}^dw_i\bf{e}_{x_i},e_G\Big):\ (w_i)_{i\leq d} \in \{0,1\}^d\Big\} \subset (\bbZ_6\wr G)_0\]
with $\{x_1,x_2,\ldots,x_d\}$ as given by
Lemmas~\ref{lem:big-equidist}, it is not hard to deduce from
Lemma~\ref{lem:approx-by-TS} how to embed copies of the Hamming cubes $\bbZ_2^d$ with uniform distortion inside $\bbZ_6\wr G$ so that their diameters grow exponentially faster than the expansion ratios of their embeddings. However, to obtain cube-\emph{quotients}
requires a little more work, in this subsection we prepare the
ground for this.

The difficulty that seems to prevent us from finding these quotients
inside $\bbZ_6\wr G$ lies in finding a single subgroup of
$\bbZ_6^{\oplus G}$ that is \emph{$G$-translation invariant} and for
which the resulting quotient will restrict to the desired quotient
on each of these Hamming cube copies separately. It is for this
reason that we make a second extension of our group, to form our
final group of the form $(\bbZ_2^{\oplus (\bbZ_6\wr G)}/V)\rtimes
(\bbZ_6\wr G)$ for a suitable choice of $(\bbZ_6\wr G)$-translation
invariant $\bbZ_2$-subspace $V$.  Even in this case we will need to
take some care in our proof that the quotient by $V$ really does
lead to the bi-Lipschitz obstructions we seek.

In this subsection we begin by putting the finite-dimensional
quotient spaces of Theorem~\ref{thm:KhoNao} into a form that is
better adapted to embedding them into a left-invariant group metric.

First, by an iterated appeal to Lemma~\ref{lem:big-equidist} we can
clearly select an increasing sequence of positive reals
$(r_n)_{n\geq 1}$, a sequence of points $(y_n)_{n\geq 1}$ in $G$ and
a sequence of pairwise-disjoint subsets $(I_n)_{n\geq 1}$ of $G$
such that
\begin{itemize}
\item $I_n \subset G \setminus \{y_1,y_2,\ldots\}$ for all $n\geq
1$,
\item $r_n \leq \rho(x,x') \leq 2Mr_n$ for all distinct $x,x' \in I_n$ and $r_n \leq \rho(x,e_G) \leq
2Mr_n$ for all $x \in I_n$, $n\geq 1$,
\item $|I_n| = d_n := 10^{r_n}$ for all $n\geq 1$, and
\item we have
\[r_1 < \rho(y_1,e_G)< r_2 < \rho(y_2,e_G) < \ldots.\]
\end{itemize}

It will matter at various points later that the sequence of positive
reals $r_1$, $\rho(y_1,e_G)$, $r_2$, $\rho(y_2,e_G)$, \ldots grows
sufficiently quickly.  However, in the remainder of this subsection
we will consider only the data of a single set $I_n$.

From this set we form the group $\bbZ_6^{I_n}$. This is a
finite-dimensional factor of $\bbZ_6^{\oplus G} \cong (\bbZ_6\wr
G)_0$.  We endow it with the Hamming metric $d_{\rm{Ham}}$ that
arises from endowing each factor copy of $\bbZ_6$ separately with
the word metric corresponding to the generating set $\{\pm 1\}$.

Within this group we will consider the subset $B_n := \{\bf{e}_x:\ x
\in I_n\}$ comprising the standard set of basis vectors. Let $C_n
\leq \bbZ_2^{I_n}$ be some $\bbZ_2$-subspace as appears in the Theorem~\ref{thm:KhoNao},
so that the cube-quotient $\bbZ_2^{I_n}/C_n^\perp$ endowed with the
quotient of the Hamming metric has $L_p$-distortion $\Omega_p(d_n)$
for all finite $p$.

The space $\bbZ_2^{I_n}$ embeds into $\bbZ_2^{\bbZ_6^{I_n}}$ through
its identification with $\bbZ_2^{B_n}$ and thence with
\[\{\W:\bbZ_6^{I_n}\to\bbZ_2:\ \spt\,\W \subseteq B_n\} \cong \bbZ_2^{B_n}\oplus
\bs{0}|_{\bbZ_6^{I_n}\setminus B_n};\] let us write
$\xi:\bbZ_2^{I_n}\into \bbZ_2^{\bbZ_6^{I_n}}$ for this embedding. In
order to avoid confusion we denote by
\[\langle \bf{a},\bf{b}\rangle := \sum_{x \in I_n}a_xb_x \in \bbZ_2\]
the $\bbZ_2$-valued inner product on $\bbZ_2^{I_n}$, and by
\[\llangle \W,\V\rrangle := \sum_{\bf{w} \in \bbZ_6^{I_n}}\W(\bf{w})\V(\bf{w}) \in \bbZ_2\]
its analog on the space $\bbZ_2^{\bbZ_6^{I_n}}$ of functions.

We will now show how to convert $C_n$ into a $\bbZ_2$-subspace $D_n
\leq \bbZ_2^{\bbZ_6^{I_n}}$ that is
\begin{itemize}
\item invariant under the translation action of $\bbZ_6^{I_n}$, and
\item such that the restriction of $\rm{TS}_{\rm{Ham}/D_n^\perp}$ to
\[\{\V + D^\perp_n:\ \V \in \bbZ_2^{B_n}\oplus \bs{0}|_{\bbZ_6^{I_n}\setminus
B_n}\}\] is equivalent under $\xi$ to the quotient of $d_{\rm{Ham}}$
on $\bbZ_2^{I_n}/C_n^\perp$.
\end{itemize}
These conclusions will be a corollary of the following lemma.

\begin{lem}\label{lem:fin-dim-subspaces}
There is a subspace $D_n \leq \bbZ_2^{\bbZ_6^{I_n}}$ that is
\begin{itemize}
\item invariant under the translation action of
$\bbZ_6^{I_n}$,
\item such that
\[D_n^\perp \cap (\bbZ_2^{B_n}\oplus \bs{0}|_{\bbZ_6^{I_n}\setminus B_n}) = \xi(C_n^\perp) \oplus \bs{0}|_{\bbZ_6^{I_n}\setminus B_n},\]
\item and such that whenever $\V,\V':\bbZ_6^{I_n}\to \bbZ_2$ are
supported on $B_n$, the distance
\[\rm{TS}_{\rm{Ham}/D_n^\perp}(\V + D_n^\perp,\V' + D_n^\perp)\]
is attained as $\rm{TS}_{\rm{Ham}}(\V,\W)$ for some $\W \in \V' +
D_n^\perp$ that is also supported on $B_n$ (and hence, by the second
point above, is identified with the sum of $\V'$ and a member of
$\xi(C_n^\perp)$).
\end{itemize}
\end{lem}

\textbf{Proof}\quad Given $\bf{w} \in \bbZ_6^T$ for any index set
$T$, let $\ol{\bf{w}} \in \bbZ_2^T$ be its image under the
coordinatewise application of the quotient homomorphism $\bbZ_6
\onto \bbZ_2$.

Now, for $\bf{a} = (a_x)_x\in \bbZ_2^{I_n}$ we define the associated
linear functional $\L_\bf{a}:\bbZ_6^{I_n}\to \bbZ_2$ by composition
with this quotient homomorphism and duality:
\[\L_\bf{a}(\bf{v}) = \langle \ol{\bf{v}},\bf{a}\rangle = \sum_{x\in I_n}\ol{v_x}a_x.\]

Given $C_n$ from Theorem~\ref{thm:KhoNao}, we simply let $D_n :=
\{\eta + \L_\bf{a}:\ \eta\in\bbZ_2,\,\bf{a} \in C_n\}$. This
collection of maps $\bbZ_6^{I_n}\to \bbZ_2$ is clearly a
$\bbZ_2$-subspace by the $\bbZ_2$-linearity of $\bf{a} \mapsto
\L_\bf{a}$, and is translation-invariant because
\[\rm{Trans}_{\bf{u}}(\L_{\bf{a}})(\bf{v}) \stackrel{\rm{dfn}}{=}
\L_{\bf{a}}(\bf{v} - \bf{u}) \stackrel{\rm{dfn}}{=} \langle
\ol{\bf{v} - \bf{u}},\bf{a}\rangle = \L_{\bf{a}}(\bf{v}) - \langle
\ol{\bf{u}},\bf{a}\rangle.\]

Observe next that if $\U,\U':\bbZ_6^{I_n}\to \bbZ_2$ are both
supported on $B_n$ then
\begin{eqnarray*}
\U - \U' \in D_n^\perp\quad&\Leftrightarrow&\quad \llangle
\U,\bs{1}_{\bbZ_6^{I_n}}\rrangle = \llangle
\U',\bs{1}_{\bbZ_6^{I_n}}\rrangle\\
&&\quad\quad\quad\hbox{and}\quad\llangle \U,\L_\bf{a}\rrangle =
\llangle
\U',\L_\bf{a}\rrangle\quad\forall \bf{a}\in C_n\\
&\Leftrightarrow&\quad \llangle \U,\L_{\bs{1}_{I_n}}\rrangle =
\llangle
\U',\L_{\bs{1}_{I_n}}\rrangle\\
&&\quad\quad\quad\hbox{and}\quad\llangle \U,\L_\bf{a}\rrangle =
\llangle
\U',\L_\bf{a}\rrangle\quad\forall \bf{a}\in C_n\\
&&\quad\quad\quad\hbox{(because $\spt\,\U,\spt\,\U'
\subseteq B_n$ and $\bs{1}_{\bbZ_6^{I_n}}|_{B_n} = \L_{\bs{1}_{I_n}}|_{B_n}$)}\\
&\Leftrightarrow&\quad\llangle \U,\L_\bf{a}\rrangle = \llangle
\U',\L_\bf{a}\rrangle\quad\forall \bf{a}\in C_n\\
&&\quad\quad\quad\hbox{(because we chose $C_n$ to contain
$\bs{1}_{I_n}$)}\\
&\Leftrightarrow&\quad \sum_{x\in I_n}\U(\bf{e}_x)\langle
\ol{\bf{e}_x},\bf{a}\rangle = \sum_{x\in
I_n}\U'(\bf{e}_x)\langle \ol{\bf{e}_x},\bf{a}\rangle\quad\forall \bf{a}\in C_n\\
&\Leftrightarrow&\quad \U - \U' \in \xi(C_n^\perp) \oplus
\bs{0}|_{\bbZ_6^{I_n}\setminus B_n}.
\end{eqnarray*}
This proves the second part of the lemma.

Now, by the addition-invariance of $\rm{TS}_{\rm{Ham}}$ on
$\bbZ_2^{\bbZ_6^{I_n}}$ (that is, invariance under addition of a
fixed member of $\bbZ_2^{\bbZ_6^{I_n}}$, \emph{not} under
translation of the base $\bbZ_6^{I_n}$), to prove the third part of
the lemma it suffices to show that for any $\V:\bbZ_6^{I_n}\to
\bbZ_2$ with $\spt\,\V \subseteq B_n$ the minimum
\[\rm{TS}_{\rm{Ham}/D_n^\perp}(D_n^\perp,\V + D_n^\perp) \stackrel{\rm{dfn}}{=} \min_{\W \in \V + D_n^\perp}\rm{TS}_{\rm{Ham}}(\bs{0},\W)\]
is attained (possibly not uniquely) for some $\W \in \V + D_n^\perp$
that also has $\spt\,\W \subseteq B_n$. Adjusting $\V$ by some
$B_n$-supported member of $D_n^\perp$ (thus, effectively, by a
member of $\xi(C_n^\perp)$) if necessary, it further suffices to
show that if $\V$ already minimizes
$\rm{TS}_{\rm{Ham}}(\bs{0},\,\cdot\,)$ among members of its
equivalence class modulo $D_n^\perp \cap (\bbZ_2^{B_n} \oplus
\bs{0}|_{\bbZ_6^{I_n}\setminus B_n})$, then any other
$\W:\bbZ_6^{I_n}\to \bbZ_2$ with $\V - \W \in D_n^\perp$ has
\[\rm{TS}_{\rm{Ham}}(\bs{0},\W) \geq \rm{TS}_{\rm{Ham}}(\bs{0},\V).\]

We will deduce this by showing how any $\W$ in $\V + D_n^\perp$ can
be explicitly adjusted to another member of this equivalence class
that is at least as close to $\bs{0}$ under $\rm{TS}_{\rm{Ham}}$ and
is supported on $B_n$.  Define the map $\R:\bbZ_2^{\bbZ_6^{I_n}}\to
\bbZ_2^{B_n}\oplus \bs{0}|_{\bbZ_6^{I_n}\setminus B_n}$ by
\[\R(\W)(\bf{e}_x) := \sum_{\bf{w} \in \bbZ_6^{I_n}}\W(\bf{w})\langle
\ol{\bf{e}_x},\ol{\bf{w}}\rangle:\] that is, $\R(\W)$ is the
$\bbZ_2$-valued indicator function of the subset of $B_n$ containing
those basis vectors that appear with odd coefficients in the
basis-decomposition of an odd number of the members of $\spt\,\W$.

Now suppose that $(\bf{x}_1 = \bs{0},\bf{x}_2,\ldots,\bf{x}_{\ell +
1} = \bs{0})$ is a cycle in $\bbZ_6^{I_n}$ that covers $\spt\,\W$.
Since $d_{\rm{Ham}}$ is a path metric on $\bbZ_6^{I_n}$
(corresponding to the usual nearest-neighbour graph on the Hamming
cube over the one-dimensional space given by $\bbZ_6$ with its
natural word metric), we may interpolate additional points along the
shortest paths joining each pair $(\bf{x}_i,\bf{x}_{i+1})$ and
re-label so that the value
\[\sum_{i=1}^\ell d_{\rm{Ham}}(\bf{x}_i,\bf{x}_{i+1})\]
in unchanged, but so that the consecutive points $\bf{x}_i$,
$\bf{x}_{i+1}$ are now neighbours in this graph.  Since the path
starts and ends at $\bs{0}$, for any $x \in I_n$ and $\bf{w}\in
\spt\,\W$ such that $w_x \neq 0$, the cycle must traverse some edge
in direction $x$ (that is, move between two points of $\bbZ_6^{I_n}$
that differ by $\pm \bf{e}_x$) at least once before reaching
$\bf{w}$, and at least once again on its way back to $\bs{0}$. It
follows that the length of the cycle is at least twice the number of
different basis vectors $\bf{e}_i$ that appear in the
basis-representation of some $\bf{w} \in \W$, and this in turn is
trivially bounded below by $|\spt\,\R(\W)|$.  On the other hand, the
same reasoning easily shows that $\rm{TS}_{\rm{Ham}}(\bs{0},\U) =
2|\spt\,\U| + 1$ whenever $\spt\,\U \subseteq B_n$ (since a cycle
that starts at $\bs{0}$ and simply goes up to each member of
$\spt\,\U$ in turn and then straight back down achieves the above
lower bound). Therefore we certainly have
\[\rm{TS}_{\rm{Ham}}(\bs{0},\W) \geq 2|\spt\,\R(\W)| + 1 = \rm{TS}_{\rm{Ham}}(\bs{0},\R(\W)).\]

It therefore suffices to show that $\R(\W) - \W \in D_n^\perp$,
since then the minimality of $\V$ among $B_n$-supported members of
its equivalence class shows that
\[\rm{TS}_{\rm{Ham}}(\bs{0},\V)\leq \rm{TS}_{\rm{Ham}}(\bs{0},\R(\W)),\]
so that concatenating this with the above inequality completes the
proof.

Thus we must show that
\[\langle \W,\eta + \L_{\bf{a}}\rangle = \langle \R(\W),\eta + \L_\bf{a}\rangle\]
for all $\eta \in \bbZ_2$ and $\bf{a} \in C_n$; in fact we will go
slightly further and show that this holds for \emph{all} $\bf{a} \in
\bbZ_2^{B_n}$.  By the $\bbZ_2$-linearity of the map $\bf{a} \mapsto
\L_\bf{a}$ it suffices to prove separately that
\[\llangle \W,\bs{1}_{\bbZ_6^{I_n}}\rrangle = \llangle \R(\W),\bs{1}_{\bbZ_6^{I_n}}\rrangle,\quad\quad\hbox{i.e.}\quad\quad|\spt\,\W| \equiv |\spt\,\R(\W)|\mod 2,\]
and that
\[\llangle \W,\L_{\ol{\bf{e}_x}}\rrangle = \llangle \R(\W),\L_{\ol{\bf{e}_x}}\rrangle\quad\quad\forall x \in I_n.\]

To prove the first of these equalities, observe that
\begin{eqnarray*}
|\spt\,\R(\W)| &=& |\{x \in I_n:\ \bf{e}_x\ \hbox{appears with odd
coeff.}\\
&&\quad\quad\quad\quad\quad\quad\quad \hbox{in an odd number
of }\bf{w}\in \spt\,\W\}|\\
&\equiv& \sum_{x \in I_n} |\{\bf{w} \in \spt\,\W:\
\langle\ol{\bf{e}}_x,\ol{\bf{w}}\rangle = 1\}| \mod 2\\
&\equiv& \sum_{x \in I_n} \sum_{\bf{w} \in \spt\,\W}\langle\ol{\bf{e}_x},\ol{\bf{w}}\rangle \mod 2\\
&\equiv& \sum_{\bf{w} \in
\bbZ_6^{I_n}}\W(\bf{w})\Big\langle \ol{\bf{w}},\sum_{x\in I_n}\ol{\bf{e}_x}\Big\rangle \mod 2\\
&\equiv& \llangle \W,\L_{\bs{1}_{I_n}}\rrangle \mod 2.
\end{eqnarray*}
Since by construction $\bs{1}_{I_n} \in C_n$ and so
$\L_{\bs{1}_{I_n}} \in D_n$, and also $\bs{1}_{\bbZ_6^{I_n}} \in
D_n$, and we know that $\W - \V \in D_n^\perp$, we have
\begin{eqnarray*}
\llangle \W,\L_{\bs{1}_{I_n}}\rrangle &=& \llangle
\V,\L_{\bs{1}_{I_n}}\rrangle\\
&\equiv& |\V| \mod 2\\
&\equiv& \llangle \V,\bs{1}_{\bbZ_6^{I_n}}\rrangle\\
&=& \llangle \W,\bs{1}_{\bbZ_6^{I_n}}\rrangle\\ &\equiv& |\spt\,\W|
\mod 2,
\end{eqnarray*}
so concatenating these equations gives the result.

To prove the second equality, we simply observe that
\begin{eqnarray*}
\llangle \W,\L_{\ol{\bf{e}_x}}\rrangle =
\sum_{\bf{w}}\W(\bf{w})\langle \ol{\bf{w}},\ol{\bf{e}_x}\rangle =
\Big\langle\sum_{\bf{w}}\sum_{y\in I_n}\W(\bf{w})\langle
\ol{\bf{e}_y},\ol{\bf{w}}\rangle\ol{\bf{e}_y},\ol{\bf{e}_x}\Big\rangle\\
= \sum_{\bf{w}}\R(\W)(\bf{w})\langle
\ol{\bf{e}_x},\ol{\bf{w}}\rangle = \llangle
\R(\W),\L_{\ol{\bf{e}_x}}\rrangle,
\end{eqnarray*}
as required. \qed

\begin{cor}\label{cor:embedding-cube-quot}
The map
\[\k^\circ_n:\bbZ_2^{I_n} \to \bbZ_2^{\bbZ_6^{I_n}}/D_n^\perp:\bf{a} \mapsto \bs{1}_{\{\bf{e}_x \in B_n:\ \langle \ol{\bf{e}}_x,\bf{a}\rangle = 1\}} + D_n^\perp\]
has kernel precisely $C_n^\perp$, and the resulting quotient map
\[\k_n:\bbZ_2^{I_n}/C_n^\perp \into \bbZ_2^{\bbZ_6^{I_n}}/D_n^\perp\]
is a bi-Lipschitz embedding of $d_{\rm{Ham}/C_n^\perp}$ into
$\rm{TS}_{\rm{Ham}/D_n^\perp}$ with distortion at most $2$ and
expansion ratio at most $2$.
\end{cor}

\textbf{Proof}\quad The identification of the kernel of $\k_n^\circ$
is already contained in the second part of
Lemma~\ref{lem:fin-dim-subspaces}, so we need only prove the
bi-Lipschitz and expansion ratio bounds.

To this end, suppose that $\bf{a},\bf{b} \in \bbZ_2^{I_n}$.  Then
the images $\k_n(\bf{a} + C_n^\perp)$ and $\k_n(\bf{b} + C_n^\perp)$
are represented modulo $D_n^\perp$ by the functions $\A :=
\bs{1}_{\{\bf{e}_x \in B_n:\ \langle \ol{\bf{e}}_x,\bf{a}\rangle =
1\}}$ and $\B := \bs{1}_{\{\bf{e}_x \in B_n:\ \langle
\ol{\bf{e}}_x,\bf{b}\rangle = 1\}}$. By the third part of
Lemma~\ref{lem:fin-dim-subspaces} the distance
\[\rm{TS}_{\rm{Ham}/D_n^\perp}(\A + D_n^\perp,\B + D_n^\perp)\]
is attained as $\rm{TS}_{\rm{Ham}}(\A',\B')$ for some $\A' \in \A +
(\xi(C_n^\perp)\oplus \bs{0}|_{\bbZ_6^{I_n}\setminus B_n})$ and $\B'
\in \B + (\xi(C_n^\perp)\oplus \bs{0}|_{\bbZ_6^{I_n}\setminus
B_n})$, and given these it is simply equal to $2|\spt(\A' + \B')| +
1$. However, we may clearly represent $\A' + D_n^\perp =
\k_n^\circ(\bf{a}')$ and $\B' + D_n^\perp = \k_n^\circ(\bf{b}')$,
and hence deduce from the first part of the corollary that $\bf{a}'
\in \bf{a} + C_n^\perp$ and $\bf{b}' \in \bf{b} + C_n^\perp$.  It
follows that $\rm{TS}_{\rm{Ham}}(\A',\B')$ is precisely
$2d_{\rm{Ham}}(\bf{a}',\bf{b}') + 1$, and thus that
\[\rm{TS}_{\rm{Ham}/D_n^\perp}(\k_n(\bf{a} + C_n^\perp),\k_n(\bf{b} + C_n^\perp)) \geq 2d_{\rm{Ham}/C_n^\perp}(\bf{a} + C_n^\perp,\bf{b} + C_n^\perp) + 1.\]
The reverse inequality also follows simply because for any $\bf{a}'
\in \bf{a} + C_n^\perp$ and $\bf{b}' \in \bf{b} + C_n^\perp$ the
images $\k_n^\circ(\bf{a}')$ and $\k_n^\circ(\bf{b}')$ are
candidates for being set equal to $\A'$ and $\B'$ above.

Thus in fact
\[\rm{TS}_{\rm{Ham}/D_n^\perp}(\k_n(\bf{a} + C_n^\perp),\k_n(\bf{b} + C_n^\perp)) = 2d_{\rm{Ham}/C_n^\perp}(\bf{a} + C_n^\perp,\bf{b} + C_n^\perp) + 1,\]
and since all distances involved are at least $1$ this gives rise to
the asserted bounds on bi-Lipschitz constant and expansion ratio.
\qed

\subsection{Completion of the proof}

We will now show how the quotients of finite-dimensional product
spaces studied in the preceding subsection can be simultaneously
bi-Lipschitzly recovered from a single translation-invariant
quotient of the metric $\rho_{((\bbZ_2^{\oplus (\bbZ_6\wr
G)})\rtimes (\bbZ_6\wr G))_0}$ on $((\bbZ_2^{\oplus (\bbZ_6\wr
G)})\rtimes (\bbZ_6\wr G))_0 = \bbZ_2^{\oplus (\bbZ_6\wr G)}$.

The tricky part is that we must find a single subspace $V \leq
\bbZ_2^{\oplus (\bbZ_6\wr G)}$ such that quotienting by it mimics
the finite-dimensional quotienting by $D_n^\perp$ studied above at
each of an increasing sequence of `scales', indexed by $n$, so that
these scales do not `interact'. This is crucial to our recovery of a
copy of
$(\bbZ_2^{\bbZ_6^{I_n}}/D_n^\perp,\rm{TS}_{\rm{Ham}/D_n^\perp})$
(and hence of $(\bbZ_2^{I_n}/C_n^\perp,d_{\rm{Ham}/C_n^\perp})$ by
Corollary~\ref{cor:embedding-cube-quot}) for each $n$. It is here
that we will see most clearly the usefulness of making a second
extension to construct our overall group.

Recall that given the group $G$ and its metric $\rho$, we can pick
positive reals $r_n$, subsets $I_n$ and points $y_n$ as in the
preceding section.  Clearly by passing to a subsequence if
necessary, we may always assume that the sequence
\[r_1 < \rho(y_1,e_G)< r_2 < \rho(y_2,e_G) < \ldots.\]
grows as fast as we please.  We will henceforth refer to this as the
\textbf{sequence of scales}.

We will now introduce the feature of $\bbZ_6$ (as opposed, say, to
$\bbZ_4$) that makes it suitable for our construction: if $1$ is its
generator, and if we let $\A$ be the $\bbZ_2$-valued indicator
function of the subset $\{0,1,3,4\}$, then $\A$ has the property (as
may easily be checked) that it and its distinct translates $\A(\cdot
- v)$ for $v \in \bbZ_6$ (of which there are only three, since
$\A(\cdot - 3) = \A$) \emph{have no linear combination modulo $2$
that is equal to the indicator function $\bs{1}_{\bbZ_6}$}. Let us
now fix this indicator function $\A$ for the rest of this paper.

Define the sequence of maps
\[\Q_n: \bbZ_6^{\oplus (G\setminus (I_1\cup I_2\cup\cdots\cup I_n\cup\{y_1,y_2,\ldots,y_n\}))}\times \bbZ_2^{\bbZ_6^{I_n}}\times \bbZ_6\to \bbZ_2^{\oplus \bbZ_6^{\oplus G}}\]
by setting
\[\Q_n(\bf{u},\W,u):\bf{w} \mapsto \delta_{\bf{u}}(\bf{w}|_{G\setminus (I_1\cup I_2\cup\cdots\cup I_n\cup\{y_1,y_2,\ldots,y_n\})})\cdot \W(\bf{w}|_{I_n})\cdot \A(w_{y_n}-u),\]
and define also
\[\t{\Q}_n: \bbZ_6^{\oplus (G\setminus (I_1\cup I_2\cup\cdots\cup I_n\cup\{y_1,y_2,\ldots,y_n\}))}\times \bbZ_2^{\bbZ_6^{I_n}}\times \bbZ_6\to \bbZ_2^{\oplus (\bbZ_6\wr G)}\]
by
\[\t{\Q}_n(\bf{u},\W,u)(\bf{w},g) = \left\{\begin{array}{ll}\Q_n(\bf{u},\W,u)(\bf{w})&\quad\hbox{if }g = e_G\\ 0&\quad\hbox{else.}\end{array}\right.\]

Thus $\Q_n(\bf{u},\W,u)$ is the $\bbZ_2$-valued indicator function
of the set of those $\bf{w} \in \bbZ_6^{\oplus G}$ such that
\begin{itemize}
\item the restriction of $\bf{w}$ to the indices in $G\setminus (I_1\cup I_2\cup\cdots\cup
I_n\cup\{y_1,y_2,\ldots,y_n\})$ agrees with $\bf{u}$,
\item the restriction $\bf{w}|_{I_n}$ is a member of $\spt\,\W$, and
\item the single coordinate $w_{y_n}$ ($= \bf{w}|_{\{y_n\}}$) is in
the translated set $\rm{Trans}_u(\spt\,\A)$.
\end{itemize}

These maps $\Q_n$ and $\t{\Q}_n$ will both underpin the construction
of $V$ and form the building blocks of our final embeddings into
$\bbZ_2^{\oplus (\bbZ_6\wr G)}$. The slightly mysterious inclusion
of a copy of $\A$ in this definition will be important for keeping
track of certain possible cancelations when we add translated copies
of images under these maps $\Q_n$, and will be clarified shortly.

Some of the immediate properties of $\Q_n$ seem worth making
explicit:

\begin{dfn}[Linearity; translation-covariance]
Each map $\Q_n$ is linear in its second argument:
\[\Q_n(\bf{u},\W,u) + \Q_n(\bf{u},\W',u) = \Q_n(\bf{u},\W + \W',u)\quad\quad\forall \W,\W' \in \bbZ_2^{\bbZ_6^{I_n}}\]
for fixed $\bf{u}$ and $u$.  We will refer to this as the
\textbf{linearity property of $\Q_n$}.

Each map $\Q_n$ also behaves well under translations in all three
arguments:
\begin{multline*}
\rm{Trans}_{\bf{w}}\Q_n(\bf{u},\W,u)(\,\cdot\,)\\ = \Q_n\big(\bf{u}
+ \bf{w}|_{G\setminus(I_1\cup\cdots\cup
I_n\cup\{y_1,\ldots,y_n\})},\rm{Trans}_{\bf{w}|_{I_n}}\W,u +
w_{y_n}\big)(\,\cdot\,)\\ \forall \bf{w},\bf{u},\W,u.
\end{multline*} We will refer to this as the \textbf{translation-covariance
property of $\Q_n$}.
\end{dfn}

Having made these definitions, we can state the last main result on
the route to Theorem~\ref{thm:main2}.

\begin{prop}\label{prop:TS-cubes}
If $(D_n)_{n\geq 1}$ is as in Lemma~\ref{lem:fin-dim-subspaces} and
the sequence of scales grows sufficiently fast then there is a
$\bbZ_2$-subspace $V \leq \bbZ_2^{\oplus (\bbZ_6\wr G)}$ that is
invariant under the coordinate right-translation action of
$\bbZ_6\wr G$ and for which the following holds: if for each $n\geq
1$ we define the mapping
\[\psi^\circ_n:\bbZ_2^{\bbZ_6^{I_n}}\to \bbZ_2^{\oplus (\bbZ_6\wr
G)}/V:\W \mapsto \t{\Q}_n(\bs{0}_{G\setminus (I_1 \cup \cdots \cup
I_n \cup \{y_1,\ldots,y_n\})},\W,0) + V,\] then
\begin{enumerate}
\item $\psi^\circ_n$ is $\bbZ_2$-linear with kernel $D_n^\perp$;
\item the resulting quotient map
\[\psi_n:\bbZ_2^{\bbZ_6^{I_n}}/D_n^\perp\into \bbZ_2^{\oplus (\bbZ_6\wr G)}/V \stackrel{\rm{inclusion}}{\into} (\bbZ_2^{\oplus (\bbZ_6\wr G)}/V)\rtimes (\bbZ_6\wr G)\]
is a bi-Lipschitz embedding of $\rm{TS}_{\rm{Ham}/D_n^\perp}$ into
$\rho_{(\bbZ_2^{\oplus (\bbZ_6\wr G)}/V)\rtimes (\bbZ_6\wr G)}$ with
distortion at most $\rm{O}(M)$ and expansion ratio $\rm{O}(r_n)$.
\end{enumerate}
\end{prop}

We will prove this proposition in several steps in this subsection,
but let us first see how it and
Corollary~\ref{cor:embedding-cube-quot} together imply
Theorem~\ref{thm:main2}, and so complete the proof of
Theorem~\ref{thm:main}.

\textbf{Proof of Theorem~\ref{thm:main2} from
Proposition~\ref{prop:TS-cubes}}\quad This follows simply by setting
\[\varphi_n := \psi_n\circ\k_n:\bbZ_2^{I_n}/C_n^\perp \into (\bbZ_2^{\oplus (\bbZ_6\wr G)}/V)\rtimes (\bbZ_6\wr G),\]
since the estimates of Corollary~\ref{cor:embedding-cube-quot} and
Proposition~\ref{prop:TS-cubes} show that $\varphi_n$ has distortion
at most $2\cdot (4M) = 8M$, which does not depend on $n$, and
expansion ratio $\rm{O}(r_n)$, where the dimension of the cube
$\bbZ_2^{I_n}$ is $d_n \geq 10^{r_n}$ and so the diameter of its
quotient $\bbZ_2^{I_n}/C_n^\perp$, which has $L_p$-distortion
proportional to this diameter, is also $\Omega(d_n)$, which grows
faster than any power of $r_n$. \qed

We will construct the subspace $V$ that we need in a number of
steps.  First, given the subspaces $D_n\leq \bbZ_2^{\bbZ_6^{I_n}}$
obtained in Lemma~\ref{lem:fin-dim-subspaces}, let
\[U_n :=
\rm{span}_{\bbZ_2}\big\{\Q_n(\bf{u},\W,u):\ \bf{u} \in
\bbZ_6^{\oplus (G\setminus (I_1 \cup \cdots \cup I_n \cup
\{y_1,\ldots,y_n\}))},\,\W \in D_n^\perp,\,u\in\bbZ_6\big\}.\]

Let us also define
\[U_n^+ :=
\rm{span}_{\bbZ_2}\big\{\Q_n(\bf{u},\W,u):\ \bf{u} \in
\bbZ_6^{\oplus (G\setminus(I_1 \cup \cdots \cup I_n \cup
\{y_1,\ldots,y_n\}))},\,\W \in
\bbZ_2^{\bbZ_6^{I_n}},\,u\in\bbZ_6\big\},\] so that $U_n \leq
U_n^+$.

It follows at once from the translation-covariance of $Q_n$ and the
translation-invariance of $D_n^\perp$ that $U_n$ (and similarly
$U_n^+$) is a translation-invariant subspace of $\bbZ_2^{\oplus
\bbZ_6^{\oplus G}}$. Now we let $U := \sum_{n\geq 1} U_n$, so that
this is still a translation-invariant $\bbZ_2$-subspace of
$\bbZ_2^{\oplus \bbZ_6^{\oplus G}}$, and finally we extend this to
the subspace
\begin{multline*}
V := \big\{\V \in \bbZ_2^{\oplus (\bbZ_6\wr G)}:\ \hbox{the map}\ (w_g)_{g\in G} \mapsto \V((w_{gg_1^{-1}})_{g\in G},g_1)\\
\hbox{lies in}\ U\ \hbox{for every}\ g_1 \in G\big\}
\end{multline*}
of $\bbZ_2^{\oplus (\bbZ_6\wr G)}$, which is manifestly invariant
under all right-translations by members of $\bbZ_6\wr G$. This is
the subspace that we will show enjoys the properties listed in
Proposition~\ref{prop:TS-cubes}.

Our first step is the following simple strengthening of
Lemma~\ref{lem:approx-by-TS}, which will reduce our proof of
Proposition~\ref{prop:TS-cubes} to a study of a simpler map into the
space $\bbZ_2^{\oplus (\bbZ_6\wr G)_0}/U$.

\begin{lem}\label{lem:first-reduction-of-embedding}
The composition of mappings
\[\begin{array}{lclcl}\bbZ_2^{\oplus (\bbZ_6\wr G)_0}/U &\into&
\bbZ_2^{\oplus (\bbZ_6\wr
G)}/V &\into& (\bbZ_2^{\oplus (\bbZ_6\wr G)}/V)\rtimes (\bbZ_6\wr G)\\
\U + U &\mapsto& (\U \oplus \bs{0}|_{(\bbZ_6\wr G)\setminus
(\bbZ_6\wr G)_0}) + V &\mapsto& \big((\U \oplus \bs{0}|_{(\bbZ_6\wr
G)\setminus (\bbZ_6\wr G)_0}) + V,(\bs{0},e_G)\big)\end{array}\] is
a well-defined $\bbZ_2$-linear injection, and moreover is a
$2$-bi-Lipschitz embedding of the metric $\rm{TS}_{\rho_{(\bbZ_6\wr
G)_0},(\bs{0},e_G)/U}$ on $\bbZ_2^{\oplus (\bbZ_6\wr G)_0}/U$ into
$\rho_{(\bbZ_2^{\oplus (\bbZ_6\wr G)}/V)\rtimes (\bbZ_6\wr G)}$ with
expansion ratio $\Theta(1)$.
\end{lem}

\textbf{Proof}\quad We have already seen in
Lemma~\ref{lem:approx-by-TS} that the second part of our
composition,
\[\begin{array}{lcl}
\bbZ_2^{\oplus (\bbZ_6\wr
G)}/V &\into& (\bbZ_2^{\oplus (\bbZ_6\wr G)}/V)\rtimes (\bbZ_6\wr G)\\
\V + V &\mapsto& \big(\V + V,(\bs{0},e_G)\big)\end{array},\] is a
$2$-bi-Lipschitz embedding of $\rm{TS}_{\rho_{(\bbZ_6\wr
G)},(\bs{0},e_G)/V}$ into $\rho_{(\bbZ_2^{\oplus (\bbZ_6\wr
G)}/V)\rtimes (\bbZ_6\wr G)}$ with expansion ratio $\Theta(1)$;
therefore it will suffice to show that the first part,
\[\bbZ_2^{\oplus (\bbZ_6\wr G)_0}/U \into
\bbZ_2^{\oplus (\bbZ_6\wr G)}/V,\] is a well-defined $\bbZ_2$-linear
injection that is an isometric embedding of
$\rm{TS}_{\rho_{(\bbZ_6\wr G)_0},(\bs{0},e_G)/U}$ into
$\rm{TS}_{\rho_{(\bbZ_6\wr G)},(\bs{0},e_G)/V}$.

Linearity is immediate, and the correctness of the definition and
injectivity hold because from the definition of $V$ we see that if
$\U, \U' \in \bbZ_2^{\oplus (\bbZ_6\wr G)}$ are both supported on
$(\bbZ_6\wr G)_0$, then they differ by a member of $V$ if and only
if $\U|_{(\bbZ_6\wr G)_0}$ and $\U'|_{(\bbZ_6\wr G)_0}$ differ by a
member of $U$.

Isometricity follows from some simple consideration of the
definition of the pinned traveling salesman metric.  Let us identify
$\bbZ_2$-valued functions on $(\bbZ_6\wr G)_0$ with their extensions
by $0$ to $\bbZ_6\wr G$ to lighten notation.  Given any $\W \in V$,
let $\W_0$ be the function $\W\cdot \bs{1}_{(\bbZ_6\wr G)_0}$
obtained by replacing $\W$ with the identically-zero function
outside the zero section $(\bbZ_6\wr G)_0$, and observe from the
definition of $V$ that $\W_0 \in V$ also.  Now if $\U$ and $\U'$ are
supported on the zero section, it follows that $\spt(\U - \U' + \W)
\supseteq \spt(\U - \U' + \W_0)$, and now appealing to the
monotonicity property of Lemma~\ref{lem:monotone} we deduce that for
such $\U$ and $\U'$ the quotient-metric distance
\[\rm{TS}_{\rho_{(\bbZ_6\wr
G)},(\bs{0},e_G)/V}(\U + V,\U' + V)\] is attained as
$\rm{TS}_{\rho_{(\bbZ_6\wr G)},(\bs{0},e_G)}(\U,\U' + \W)$ for some
$\W$ also supported on $(\bbZ_6\wr G)_0$, and hence agrees with
\[\rm{TS}_{\rho_{(\bbZ_6\wr
G)_0},(\bs{0},e_G)/U}(\U + U,\U' + U),\] as required. \qed

We next prove the two lemmas that will underpin the main estimates
involved in the proof of Proposition~\ref{prop:TS-cubes}.

\begin{lem}\label{lem:lin-indep-+}
Suppose that $\eta \in \bbZ_2$, $m\geq 1$ is an integer and that we
are given
\begin{itemize}
\item a fixed point $\bf{u}^\circ \in \bbZ_6^{\oplus (G\setminus(I_1 \cup\cdots\cup I_m\cup\{y_1,\ldots,y_m\}))}$, and
\item for each $n \leq m$, a finite family of functions
\[\F_n := \{\Q_n(\bf{u}^{(n)}_1,\W^{(n)}_1,u^{(n)}_1),\Q_n(\bf{u}^{(n)}_2,\W^{(n)}_2,u^{(n)}_2),\ldots,\Q_n(\bf{u}^{(n)}_{i_n},\W^{(n)}_{i_n},u^{(n)}_{i_n})\}\]
from $U_n^+$ such that $\bf{u}^{(n)}_i|_{(G\setminus(I_1
\cup\cdots\cup I_m\cup\{y_1,\ldots,y_m\}))} = \bf{u}^\circ$ for all
$i \leq i_n$
\end{itemize}
such that
\[\sum_{n=1}^m\sum_{i=1}^{i_n}\Q_n(\bf{u}^{(n)}_i,\W^{(n)}_i,u^{(n)}_i) = \eta \bs{1}_{\{\bf{u}^\circ\} \times \bbZ_6^{I_1 \cup \cdots \cup I_m \cup \{y_1,\ldots,y_m\}}}.\]
Then we must have $\eta = 0$ and
\[\sum_{i=1}^{i_n}\Q_n(\bf{u}^{(n)}_i,\W^{(n)}_i,u^{(n)}_i) = 0\]
for each $n\leq m$ separately.
\end{lem}

\textbf{Remark}\quad It will be important that we have this lemma
available for basis members of the larger spaces $U_n^+$, not just
of $U_n$. \fin

\textbf{Proof}\quad This will follow by induction on the number
among the families $\F_1$, $\F_2$, \ldots, $\F_m$ that are nonempty.

\quad\textbf{Base clause}\quad Suppose that only one $\F_n$ is
nonempty, in which case our assumption reads
\[\sum_{i=1}^{i_n}\Q_n(\bf{u}^{(n)}_i,\W^{(n)}_i,u^{(n)}_i) = \eta \bs{1}_{\{\bf{u}^\circ\} \times \bbZ_6^{I_1 \cup \cdots \cup I_m \cup \{y_1,\ldots,y_m\}}},\]
and we need only show that $\eta = 0$.  We prove this by
contradiction, so suppose instead that $\eta = 1$.  The function
$\Q_n(\bf{u}^{(n)}_i,\W^{(n)}_i,u^{(n)}_i)$ is supported on the
(infinite-dimensional) cylinder $\{\bf{u}^{(n)}_i\} \times
\bbZ_6^{(I_1 \cup \cdots \cup I_n \cup\{y_1,\ldots,y_n\})}$, and so
if $\bf{u}^{(n)}_i \neq \bf{u}^{(n)}_j$ then
$\Q_n(\bf{u}^{(n)}_i,\W^{(n)}_i,u^{(n)}_i)$ and
$\Q_n(\bf{u}^{(n)}_j,\W^{(n)}_j,u^{(n)}_j)$ are disjointly
supported.  Now suppose that $\bf{v}$ is a vector that appears in
the list $\bf{u}^{(n)}_1$, \ldots, $\bf{u}^{(n)}_{i_n}$ (so that by
assumption, $\bf{v}$ is an extension of $\bf{u}^\circ$), and let
$j_1 < j_2 < \ldots < j_k$ be those values of $i \leq i_n$ where it
appears. It follows from the above equation that
\[\sum_{s = 1}^k\Q_n(\bf{v},\W^{(n)}_{j_s},u^{(n)}_{j_s}) =
1\quad\quad\hbox{on }\{\bf{v}\} \times \bbZ_6^{I_n \cup
\{y_n\}}\times \{\bs{0}_{I_1 \cup \cdots \cup I_{n-1}
\cup\{y_1,\ldots,y_{n-1}\}}\}.\]

Now pick any $\bf{w} \in \bbZ_6^{I_n}$; the above equation requires
that $\W^{(n)}_{j_s}(\bf{w}) = 1$ for at least one $s \leq k$.
Letting $s_1 < s_2 < \ldots < s_\ell$ be those values of $s \leq k$
where $\W^{(n)}_{j_s}(\bf{w}) = 1$, we deduce from the above
equation that
\[\sum_{r =
1}^\ell\Q_n(\bf{v},\W^{(n)}_{j_{s_r}},u^{(n)}_{j_{s_r}}) =
1\quad\quad\hbox{on }\{\bf{v}\}\times\{\bf{w}\} \times \bbZ_6\times
\{\bs{0}_{I_1 \cup \cdots \cup I_{n-1}
\cup\{y_1,\ldots,y_{n-1}\}}\}.\]

However, each of the functions
$\Q_n(\bf{v},\W^{(n)}_{j_{s_r}},u^{(n)}_{j_{s_r}})$ restricted to
the set $\{\bf{v}\oplus \bf{w}\} \times \bbZ_6\times \{\bs{0}_{I_1
\cup \cdots \cup I_{n-1} \cup\{y_1,\ldots,y_{n-1}\}}\} \cong \bbZ_6$
is simply a rotated copy of the function $\A$, and we introduced
this precisely so as to have the property that the constant function
$1$ cannot be made as a linear combination of its rotates.  This
gives the desired contradiction.

\quad\textbf{Recursion clause}\quad Suppose that $1 \leq n_1 < n_2 <
\ldots < n_\ell = m$ with $\ell \geq 2$ are the values of $n$ for
which $\F_n$ is nonempty.  We will show that in this case the
assumption that
\[\sum_{j=1}^\ell\sum_{i=1}^{i_{n_j}}\Q_{n_j}(\bf{u}^{(n_j)}_i,\W^{(n_j)}_i,u^{(n_j)}_i) = \eta \bs{1}_{\{\bf{u}^\circ\} \times \bbZ_6^{I_1 \cup \cdots \cup I_m \cup \{y_1,\ldots,y_m\}}},\]
implies that also
\[\sum_{j=2}^\ell\sum_{i=1}^{i_{n_j}}\Q_{n_j}(\bf{u}^{(n_j)}_i,\W^{(n_j)}_i,u^{(n_j)}_i) = \eta \bs{1}_{\{\bf{u}^\circ\} \times \bbZ_6^{I_1 \cup \cdots \cup I_m \cup \{y_1,\ldots,y_m\}}},\]
(that is, omitting the first term of the outer sum) so that an
induction on $\ell$ completes the proof.

To see this, observe that each of the functions
$\Q_{n_1}(\bf{u}^{(n_1)}_i,\W^{(n_1)}_i,u^{(n_1)}_i)$ is supported
on the set $\{\bf{u}_i^{(n_1)}\}\times \bbZ_6^{I_{n_1} \cup I_{n_1 -
1}\cup\cdots \cup I_1 \cup \{y_{n_1},y_{n_1-1},\ldots,y_1\}}$,
whereas each of the functions
$\Q_{n_j}(\bf{u}^{(n_j)}_i,\W^{(n_j)}_i,u^{(n_j)}_i)$ for $j\geq 2$
is constant on every set of the form $\{\bf{v}\}\times
\bbZ_6^{I_{n_1}\cup\cdots \cup I_1 \cup \{y_{n_1},\ldots,y_1\}}$
with $\bf{v} \in \bbZ_6^{\oplus G\setminus(I_{n_1}\cup\cdots \cup
I_1 \cup \{y_{n_1},\ldots,y_1\})}$, as is the function
$\eta\bs{1}_{\{\bf{u}^\circ\} \times \bbZ_6^{I_1 \cup \cdots \cup
I_m \cup \{y_1,\ldots,y_m\}}}$.  Therefore the stated cancelation
can take place only if for every $\bf{v}$ that appears in the list
$\bf{u}^{(n_1)}_1$, $\bf{u}^{(n_1)}_2$, \ldots,
$\bf{u}^{(n_1)}_{i_{n_1}}$ we have
\[\sum_{i \leq i_{n_1},\,\bf{u}^{(n_1)}_i = \bf{v}}\Q_{n_1}(\bf{v},\W^{(n_1)}_i,u^{(n_1)}_i) = \rm{const.}\quad\quad\hbox{on }\{\bf{v}\}\times \bbZ_6^{I_{n_1}\cup\cdots \cup I_1 \cup
\{y_{n_1},\ldots,y_1\}},\] and applying the argument for the base
clause to this sub-sum tells us that the constant in question must
be zero.  Summing over $\bf{v}$ now gives
\[\sum_{i = 1}^{i_{n_1}}\Q_{n_1}(\bf{u}^{(n_1)}_i,\W^{(n_1)}_i,u^{(n_1)}_i) = 0,\]
and finally subtracting this sum from the initially-given equation
completes the recursion step, so the induction continues. \qed

\textbf{Remark}\quad If we assume that $\eta = 0$ a priori then the
above proposition asserts simply that the subspaces $U_n^+ \leq
\bbZ_2^{\oplus \bbZ_6^{\oplus G}}$, $n=1,2,\ldots$, are linearly
independent over $\bbZ_2$.  Indeed, given any equation
\[\sum_{n=1}^m\sum_{i=1}^{i_n}\Q_n(\bf{u}^{(n)}_i,\W^{(n)}_i,u^{(n)}_i) = 0\]
we may first partition the left-hand side into sub-sums according to the values of $\bf{u}^{(n)}_i|_{(G\setminus(I_1
\cup\cdots\cup I_m\cup\{y_1,\ldots,y_m\}))}$, and now these sub-sums must all vanish separately because any two maps $\Q_n(\bf{u}^{(n)}_i,\W^{(n)}_i,u^{(n)}_i)$ and $\Q_{n'}(\bf{u}^{(n')}_j,\W^{(n')}_j,u^{(n')}_j)$ have disjoint support if these restrictions do not agree.  Each of these sub-sums falls within the conditions of the proposition.  \fin

We are now ready to deduce the crucial relation between the metric
$\rm{TS}_{\rho_{(\bbZ_6\wr G)_0},(\bs{0},e_G)}$ on $\bbZ_2^{\oplus
(\bbZ_6\wr G)_0}$ and the metrics $\rm{TS}_{\rm{Ham}}$ on some of
its finite-dimensional factor spaces.

\begin{lem}\label{lem:ignore-cluster-size}
Provided the sequence of scales grows sufficiently fast, the metric
$\rm{TS}_{\rho_{(\bbZ_6\wr G)_0},(\bs{0},e_G)}$ admits the following
kind of approximation. Suppose that
$\U_n,\V_n:\bbZ_6^{I_n}\to\bbZ_2$ are non-zero and distinct, and
that $\U,\V:\bbZ_6^{\oplus G}\to \bbZ_2$ are such that
\[\spt\,\U \subseteq \{\bs{0}\}\times \spt\,\U_n\times \bbZ_6^{(I_1\cup\cdots\cup I_{n-1}\cup\{y_1,\ldots,y_n\})},\]
\[\spt\,\V \subseteq \{\bs{0}\}\times \spt\,\V_n\times \bbZ_6^{(I_1\cup\cdots\cup I_{n-1}\cup\{y_1,\ldots,y_n\})},\]
and
\[\spt\,\U\cap \big(\{\bs{0}\}\times \{\bf{u}\}\times \bbZ_6^{(I_1\cup\cdots\cup I_{n-1}\cup\{y_1,\ldots,y_n\})}\big) \neq \emptyset\quad\quad\forall \bf{u} \in \spt\,\U_n,\]
\[\spt\,\V\cap \big(\{\bs{0}\}\times \{\bf{u}\}\times \bbZ_6^{(I_1\cup\cdots\cup I_{n-1}\cup\{y_1,\ldots,y_n\})}\big) \neq \emptyset\quad\quad\forall \bf{u} \in \spt\,\V_n,\]
\begin{multline*}
\U|_{\{\bs{0}\}\times \{\bf{u}\}\times \bbZ_6^{(I_1\cup\cdots\cup
I_{n-1}\cup\{y_1,\ldots,y_n\})}} = \V|_{\{\bs{0}\}\times
\{\bf{u}\}\times \bbZ_6^{(I_1\cup\cdots\cup
I_{n-1}\cup\{y_1,\ldots,y_n\})}}\\ \forall \bf{u} \in \spt\,\U_n\cap
\spt\,\V_n.
\end{multline*}

Then we have
\[r_n\cdot \rm{TS}_{\rm{Ham}}(\U_n,\V_n) \leq \rm{TS}_{\rho_{(\bbZ_6\wr G)_0},(\bs{0},e_G)}(\U,\V) \leq 4Mr_n\cdot \rm{TS}_{\rm{Ham}}(\U_n,\V_n)\]
(irrespective of the actual sizes of the sets $\spt\,\U$,
$\spt\,\V$).
\end{lem}

\textbf{Proof}\quad We now write simply $\bs{0}$ for the zero of
either $\bbZ_6^{\oplus G}$ or $\bbZ_6^{\oplus G\setminus
(I_1\cup\cdots\cup I_n \cup\{y_1,\ldots,y_n\})}$ for any $n$, since
among these possibilities the index set will always be clear from
the context.

If the distance $r_n$ is sufficiently large compared with $r_1$,
$\rho(e_G,y_1)$, $r_2$, $\rho(e_G,y_2)$, \ldots, $r_{n-1}$ and
$\rho(e_G,y_n)$, and if $\bf{u},\bf{u}' \in \bbZ_6^{I_n}$ are
distinct, then the pinned traveling salesman distance in $(\bbZ_6\wr
G)_0$ between any points of the cylinder sets
\[\{\bs{0}\}\times \{\bf{u}\}\times \bbZ_6^{(I_1\cup\cdots\cup I_{n-1}\cup\{y_1,\ldots,y_n\})}\]
and
\[\{\bs{0}\}\times \{\bf{u}'\}\times \bbZ_6^{(I_1\cup\cdots\cup I_{n-1}\cup\{y_1,\ldots,y_n\})}\]
will be
\[c\cdot \rho_{(\bbZ_6\wr G)_0}(\bs{0} \oplus \bf{u} \oplus \bs{0}_{I_1\cup\cdots\cup I_{n-1}\cup\{y_1,\ldots,y_n\}},\bs{0} \oplus \bf{u}' \oplus \bs{0}_{I_1\cup\cdots\cup I_{n-1}\cup\{y_1,\ldots,y_n\}}) \geq r_n,\]
for some $1 \leq c \leq 2$, since the length of a cycle in $G$
needed to cover all the coordinates where $\bf{u}$ differs from
$\bf{u}'$ dwarfs the maximum number of steps that could possibly be
needed to cover all the differences between two points necessary at
coordinates indexed by $I_1\cup\cdots\cup
I_{n-1}\cup\{y_1,\ldots,y_n\}$.

Similarly, if $r_n$ is sufficiently large then the pinned traveling
salesman distance between these cylinder sets also dwarfs the
maximum pinned traveling salesman length of any
non-self-intersecting path within either of these cylinder sets.

It follows for $\U$ and $\V$ as given that any cycle in $(\bbZ_6\wr
G)_0$ starting and ending at $(\bs{0},e_G)$ and covering $\spt(\U +
\V)$ needs to cover at least one point of each of the sets
\[\spt\,\U\cap
\big(\{\bs{0}\}\times \{\bf{u}\}\times \bbZ_6^{(I_1\cup\cdots\cup
I_{n-1}\cup\{y_1,\ldots,y_n\})}\big) \quad\quad\hbox{for }\bf{u} \in
\spt\,\U_n\setminus \spt\,\V_n\] and
\[\spt\,\V\cap
\big(\{\bs{0}\}\times \{\bf{u}\}\times \bbZ_6^{(I_1\cup\cdots\cup
I_{n-1}\cup\{y_1,\ldots,y_n\})}\big) \quad\quad\hbox{for }\bf{u} \in
\spt\,\V_n\setminus \spt\,\U_n,\] but no others, and that the
distance between any two of these sets is much larger than the
maximum number of steps in a path that could be needed to cover the
necessary points within any one of them.

Therefore the length of such a covering traveling salesman
cycle for $\spt(\U + \V)$ is bounded below by the length of a
traveling salesman cycle in $(\bbZ_6\wr G)_0$ for the set
\[\big\{\bs{0} \oplus \bf{u}\oplus \bs{0}|_{I_1\cup\cdots\cup
I_{n-1}\cup\{y_1,\ldots,y_n\}}:\ \bf{u} \in \spt\,\U_n\triangle
\spt\,\V_n\big\},\] and if $r_n$ is sufficiently large compared with
its predecessors then it may also be bounded from above by twice
this number.  Finally, since by construction any two index points of
$I_n$ are separated by a $\rho$-distance that lies in $[r_n,2Mr_n]$,
and so any covering cycle for this latter set must traverse an
additional distance of at least $r_n$ and at most $2Mr_n$ for each
new point that it must visit, we deduce this latter distance
lies between $r_n\cdot \rm{TS}_{\rm{Ham}}(\U_n,\V_n)$ and
$2Mr_n\cdot \rm{TS}_{\rm{Ham}}(\U_n,\V_n)$.

Combining the bounds obtained above now completes the proof. \qed

\textbf{Proof of Proposition~\ref{prop:TS-cubes}}\quad First,
Lemma~\ref{lem:first-reduction-of-embedding} allows us to consider
instead the maps
\[\l^\circ_n:\bbZ_2^{\bbZ_6^{I_n}}\to \bbZ_2^{\oplus (\bbZ_6\wr G)_0}/U: \W \mapsto \Q_n(\bs{0},\W,0) + U,\]
for which it will suffice to show the corresponding properties:
\begin{enumerate}
\item $\l_n^\circ$ is $\bbZ_2$-linear with kernel $D_n^\perp$;
\item the resulting quotient map
\[\l_n:\bbZ_2^{\bbZ_6^{I_n}}/D_n^\perp\into \bbZ_2^{\oplus (\bbZ_6\wr G)_0}/U\]
is a bi-Lipschitz embedding of $\rm{TS}_{\rm{Ham}/D_n^\perp}$ into
$\rm{TS}_{\rho_{(\bbZ_6\wr G)_0},(\bs{0},e_G)/U}$ with distortion at
most $\rm{O}(M)$ and expansion ratio $\rm{O}(r_n)$.
\end{enumerate}

\quad\textbf{1.}\quad It is clear from the definition of $U_n \leq
U$ that $\l^\circ_n$ annihilates $D_n^\perp$, so we need only show
that it does not annihilate any larger subspace of
$\bbZ_2^{\bbZ_6^{I_n}}$. However, a special case of
Lemma~\ref{lem:lin-indep-+} tells us that the subspaces $U^+_n \leq
\bbZ_2^{\oplus (\bbZ_6\wr G)_0}$ are linearly independent (as
remarked immediately after that lemma). It follows that if
$\l^\circ_n(\V) = \Q_n(\bs{0},\V,0) + U$ is equal to the identity
element $U$ in $\bbZ_2^{\oplus (\bbZ_6\wr G)_0}/U$, then there are
elements $\A_1 \in U_1$, $\A_2 \in U_2$, \ldots, $\A_m \in U_m$ for
some $m \geq n$ such that
\[\A_1 + \A_2 + \cdots + \big(\Q_n(\bs{0},\V,0) + \A_n\big) + \cdots + \A_m =
\bs{0},\] and this is possible only if each $\A_i$ is individually
zero for $i \neq n$ and also $\Q_n(\bs{0},\V,0) = \A_n$.

We now express $\A_n$ as
\[\sum_{i=1}^{i_1} \Q_n(\bf{u}_i,\W_i,u_i)\quad\quad \W_i \in D_n^\perp,\,u_i \in \bbZ_6\quad\forall i \leq i_1,\]
and will show that we can have $\Q_n(\bs{0},\V,0)$ equal to such a
sum only if in fact $\V$ is itself a member of $D_n^\perp$. First,
we may clearly discard all terms of this sum for which $\bf{u}_i
\neq \bs{0}$, since Lemma~\ref{lem:lin-indep-+} tells us immediately
that these must cancel to $0$, and relabel so that $\bf{u}_i =
\bs{0}$ for every $i\leq i_1$. Therefore, recalling the definition
of $\Q_n$ and omitting the fixed vector $\bs{0}$, the above equation
becomes
\[\V(\bf{w})\cdot \A(w) = \sum_{i=1}^{i_1} \W_i(\bf{w})\cdot\A(w - u_i)\quad\quad\forall \bf{w} \in \bbZ_6^{I_n},\, w \in \bbZ_6.\]

However, if we know simply fix $w \in \spt\,\A$, this equation
simplifies to
\[\V(\bf{w}) = \sum_{1 \leq i \leq i_1,\,\A(w - u_i) = 1} \W_i(\bf{w})\quad\quad\forall \bf{w} \in \bbZ_6^{I_n},\]
and so we have expressed $\V$ as a linear combination of members of
$D_n^\perp$, and hence proved that it itself lies in $D_n^\perp$, as
required.

\quad\textbf{2.}\quad We prove the two necessary inequalities
separately. First we prove that
\[\rm{TS}_{\rho_{(\bbZ_6\wr
G)_0},(\bs{0},e_G)/U}(\l_n(\U + D_n^\perp),\l_n(\V+D_n^\perp)) \leq
4Mr_n\cdot \rm{TS}_{\rm{Ham}/D_n^\perp}(\U + D_n^\perp,\V +
D_n^\perp).\] Indeed, if $\W \in D_n^\perp$ minimizes
$\rm{TS}_{\rm{Ham}}(\U,\V + \W)$, then $\Q_n(\bs{0},\W,0) \in U_n
\leq U$ and by a simple application of
Lemma~\ref{lem:ignore-cluster-size} the distance
\[\rm{TS}_{\rho_{(\bbZ_6\wr
G)_0},(\bs{0},e_G)}(\Q_n(\bs{0},\U,0),\Q_n(\bs{0},\V,0) +
\Q_n(\bs{0},\W,0))\] is at least $r_n$ and at most $4Mr_n$ times
$\rm{TS}_{\rm{Ham}}(\U,\V + \W)$, as required.

The necessary reverse inequality
\[\rm{TS}_{\rho_{(\bbZ_6\wr
G)_0},(\bs{0},e_G)/U}(\l_n(\U + D_n^\perp),\l_n(\V + D_n^\perp))
\gtrsim r_n\cdot \rm{TS}_{\rm{Ham}/D_n^\perp}(\U + D_n^\perp,\V +
D_n^\perp)\] requires a little more work.  By addition-invariance we
may assume that $\U = \bs{0}$.  Now suppose that $\V \in
\bbZ_2^{\bbZ_6^{I_n}}$ and $m\geq n$ and that
\begin{multline*}
\F_k :=
\{\Q_k(\bf{u}^{(k)}_1,\W^{(k)}_1,u^{(k)}_1),\Q_k(\bf{u}^{(k)}_2,\W^{(k)}_2,u^{(k)}_2),\ldots,
\Q_k(\bf{u}^{(k)}_{i_k},\W^{(k)}_{i_k},u^{(k)}_{i_k})\}\\ \hbox{for
}k=1,2,\ldots,m
\end{multline*}
are finite collections in each $U_k$.  Clearly we may assume in our notation that the functions
$\W^{(k)}_i$ are all non-zero, that $\F_m \neq \emptyset$ and
finally that $m$ is minimal subject to all these other restrictions.

To complete the proof we will show that for some such finite data the combined
collection $\F := \bigcup_{k\leq m}\F_k$ actually minimizes the value
\begin{eqnarray}\label{eq:minimized}
&&\rm{TS}_{\rho_{(\bbZ_6\wr
G)_0},(\bs{0},e_G)}\big(\bs{0},\Q_n(\bs{0},\V,0) + \S\F\big)
\end{eqnarray}
over all such finite tuples of finite collections, where we
write $\S\F$ simply for the sum in $\bbZ_2^{\oplus (\bbZ_6\wr G)_0}$
of all elements of $\F$, and then from this minimizing $\F$ we will
construct some $\W \in D_n^\perp$ so that
\[\rm{TS}_{\rho_{(\bbZ_6\wr G)_0},(\bs{0},e_G)}\big(\bs{0},\Q_n(\bs{0},\V,0) + \S\F\big) \geq r_n\cdot \rm{TS}_{\rm{Ham}}(\bs{0},\V + \W).\]

\quad\textbf{Step (i)}\quad We first show by contradiction that in order to infimize the expression~(\ref{eq:minimized}) over finite tuples $\F_k$, $k\leq m$, it suffices to assume that $m \leq
n$. Indeed, suppose instead that $m\geq n+1$. First, if any $\bf{u}^{(m)}_i$
is not $\bs{0}$, then we may simply omit the corresponding member of
$\F_m$ together with all members of any $\F_k$ with $k < m$ whose
supports it dominates, and obtain a new family $\F'$ for which
\[\rm{TS}_{\rho_{(\bbZ_6\wr G)_0},(\bs{0},e_G)}\big(\bs{0},\Q_n(\bs{0},\V,0) + \S\F'\big) < \rm{TS}_{\rho_{(\bbZ_6\wr G)_0},(\bs{0},e_G)}\big(\bs{0},\Q_n(\bs{0},\V,0) + \S\F\big).\]

Therefore it suffices to assume that $\bf{u}^{(m)}_i = \bs{0}$ for all $i \leq i_m$. Next, for
any given $w \in \bbZ_6$, we have
\begin{multline*}
\sum_{i \leq i_m}\Q_m(\bs{0},\W^{(m)}_i,u^{(m)}_i)(\bs{0}\oplus
\bf{w}\oplus w\oplus \bf{u}) = \sum_{i\leq i_m,\,\A(w - u^{(m)}_i) =
1}\W^{(m)}_i(\bf{w})\\ \hbox{for }\bf{w} \in \bbZ_6^{I_m},\,\bf{u}
\in \bbZ_6^{I_1 \cup \cdots \cup I_{m-1} \cup
\{y_1,\ldots,y_{m-1}\}}.
\end{multline*}
Regarded as a function of $\bf{w}\in \bbZ_6^{I_m}$ alone, we see
that this is a member of $D_m^\perp$.  If it is zero for every $w
\in \bbZ_6$, then summing over $w$ shows that $\S\F_m = 0$, and so
we may simply discard $\F_m$ to leave a new family $\F' :=
\bigcup_{k \leq m-1}\F_k$ that achieves the same value for the
expression~(\ref{eq:minimized}) but has a smaller value of $m$.  Hence after making finitely many such omissions, we may assume also that there is some $w \in
\bbZ_6$ for which the above sum specifies a nonzero member of
$D_m^\perp$; let us call that member $\Y$.  Then we also have $\Y
\neq \delta_{\bs{0}}$ (simply because $\delta_{\bs{0}} \not\in
D_m^\perp$), and so to this $w$ there corresponds some $\bf{w} \in
\bbZ_6^{I_m}\setminus \{\bs{0}\}$ for which $\Y(\bf{w}) \neq 0$.

However, it now follows that the function $\sum_{i \leq
i_m}\Q_m(\bs{0},\W^{(m)}_i,u^{(m)}_i)$ is non-zero on the cylinder
$\{\bs{0}\}\times \{\bf{w}\}\times \{w\}\times
\bbZ_6^{I_1\cup\cdots\cup I_{m-1}\cup \{y_1,\ldots,y_{m-1}\}}$, and
so actually takes the constant value $1$ on that cylinder.  By
Lemma~\ref{lem:lin-indep-+} the sum function $\Q_n(\bs{0},\V,0) +
\S(\F_1\cup\cdots\cup \F_{m-1})$ cannot be constant and equal to $1$
on this cylinder, and so overall the function $\Q_n(\bs{0},\V,0) + \S\F$
does not vanish on this cylinder. On the other hand, if our sequence
of scales grows fast enough, then any point of this cylinder is by
itself much further away from $(\bs{0},e_G)$ in $\rho_{(\bbZ_6\wr
G)_0}$ than any possible value of $\rm{TS}_{\rho_{(\bbZ_6\wr
G)_0,(\bs{0},e_G)}}(\bs{0},\Q_n(\bs{0},\V,0))$, and so in this case even the empty family gives a smaller value for the expression~(\ref{eq:minimized}) than does $\F$.  Therefore the infimum of the values~(\ref{eq:minimized}) is unchanged if we only infimize over families of collections $\F$ for which $m\leq n$.

\quad\textbf{Step (ii)}\quad Having reduced to the case $m \leq n$, we
can now argue as above to reduce the evaluation of the infimum of~(\ref{eq:minimized}) further to the case when
\[\bf{u}^{(k)}_i|_{G\setminus (I_1\cup\cdots\cup I_n\cup\{y_1,\ldots,y_n\})} = \bs{0}_{G\setminus (I_1\cup\cdots\cup I_n\cup\{y_1,\ldots,y_n\})}\quad\quad\forall k\leq m,\]
since any of the functions in question is either supported on the
cylinder specified by this collection of $0$-valued coordinates or
is identically $0$ on it, and so discarding those that are
identically zero gives another collection $\F$ for which the value~(\ref{eq:minimized}) is not any larger.  This last reduction leaves only finitely many possibilities for the collection $\F$, and so now we know that we may pick one that is actually minimizing.  One last `processing' step will give us the construction of $\W$ from it.

\quad\textbf{Step (iii)}\quad For this, it now follows from another application of
Lemma~\ref{lem:lin-indep-+} that on any of the cylinders
$\{\bf{0}\}\times \{\bf{w}\}\times \{0\}\times
\bbZ_6^{I_1\cup\cdots\cup I_{n-1}\cup \{y_1,\ldots,y_{n-1}\}}$ for
$\bf{w} \in \bbZ_6^{I_n}$ the sum function $\Q_n(\bs{0},\V,0) +
\S\F$ can vanish only if each individual sum $\S\F_k$ vanishes there
for each $k \leq n-1$, and also $\Q_n(\bs{0},\V,0) + \S\F_n$
vanishes there. Substituting the definition of $\Q_n$ into this
latter condition, it becomes that
\[\V(\bf{w}) + \sum_{1 \leq i \leq i_n,\,\A(-u^{(n)}_i) = 1}\W^{(n)}_i(\bf{w}) = 0.\]

Therefore, the indicator function of the set \[\big\{\bf{w} \in
\bbZ_6^{I_n}: \big(\Q_n(\bs{0},\V,0) + \S\F\big)|_{\{\bf{0}\}\times
\{\bf{w}\}\times \{0\}\times \bbZ_6^{I_1\cup\cdots\cup I_{n-1}\cup
\{y_1,\ldots,y_{n-1}\}}} \not\equiv 0\big\}\] is precisely $\V +
\sum_{1 \leq i \leq i_n,\,\A(-u^{(n)}_i) = 1}\W^{(n)}_i$, and
therefore it is a member of $\V + D_n^\perp$; let us call it $\V +
\W$. By Lemma~\ref{lem:monotone} and
Lemma~\ref{lem:ignore-cluster-size} this now implies at once that
\[\rm{TS}_{\rho_{(\bbZ_6\wr G)_0},(\bs{0},e_G)}\big(\bs{0},\Q_n(\bs{0},\V,0) + \S\F\big) \geq r_n\cdot \rm{TS}_{\rm{Ham}}(\bs{0},\V + \W).\]
as required. \qed

\section{Further questions}

I suspect that the construction above can be extended so far as to
give a finitely-generated amenable group $G$ with a word metric
$\rho$ for which any $1$-Lipschitz embedding $f:G \into L_p$, $p\in
[1,\infty)$, must be such that
\[\|f(g) - f(h)\| \lesssim \log \rho(g,h)\]
for some collection of pairs $g,h \in G$ among which $\rho(g,h)$ can
be arbitrarily large.

\textbf{Question}\quad Does every finitely-generated amenable group $G$ admit Lipschitz embeddings $f:G\into L_p$ for every $p \in [1,\infty)$ such that $\|f(g) - f(h)\|_p \gtrsim \log \rho(g,h)$ for all $g,h \in G$?

Applied with a little less violence (that is, with $\eta < 1$ in
Lemma~\ref{lem:obs}), the methods of this paper may also be useful
for finding examples of finitely-generated amenable groups with
specified compression exponents in $(0,1)$, although of course these
methods will need to be complemented with matching lower-bound
proofs for the groups in question (that is, constructions of particular good-distortion embeddings, on which the ideas of the present paper do not seem to bear directly).  The kind of construction used above may also provide interesting test-cases for Question 10.6 of Naor and Peres~\cite{NaoPer09}, which asks for more general
methods for estimating compression exponents of semidirect products;
and on the question of which values are realizable as Euclidean
compression exponents for finitely-generated amenable groups,
specializing the result of Arzhantseva, Drutu and Sapir
in~\cite{ArzDruSap09} that all values in $[0,1]$ can be realized as
the compression exponents of finitely-generated groups that are not
necessarily amenable.

It is interesting to note that in~\cite{KhoNao06} Khot and Naor use
inequalities described in terms of random walks on Hamming cubes to
prove Theorem~\ref{thm:KhoNao}. Thus, although the random walk upper
bound on compression exponents given by Naor and Peres cannot reach
the zero-compression-exponent regime of the group $(\bbZ_2^{\oplus
H}/V)\rtimes H$ constructed above, it seems that a modification of
their method in which we allow ourselves to consider a sequence of
random walks on our group supported on the images of these embedded
finite cube-quotients would translate into the correct zero upper
bound on the compression exponent.  It would be interesting to find
an example of a finitely-generated amenable group for which some
other obstruction to embeddings with positive compression exponents
is needed, genuinely unrelated to inequalities concerning random walks.

One candidate for such an obstruction would be a sequence of
embedded copies of the cubes $(\bbZ_{2m}^n,\ell_\infty)$: given an
embedded sequence of these with expansion ratios not growing too
fast, we could instead use an analog of Lemma~\ref{lem:obs} based on the notion of non-linear cotype introduced by Mendel and Naor in~\cite{MenNao08}, which is
characterized in terms of these finite spaces.  Given a suitable
tradeoff between the sizes and expansion ratios, the arguments
leading to Mendel and Naor's Theorem 1.11 of~\cite{MenNao08} would
imply that such a sequence of embedded $\ell_\infty$-cubes obstructs
good-compression embeddings into any Banach space with nontrivial
type and cotype $q < \infty$, provided in addition that the
parameter $m$ of the embedded cubes can be chosen to grow faster
than $n^{1/q}$.  (Note also that whether the assumption of nontrivial
type is necessary is one of the major outstanding problems from
their paper.) However, at this stage I do not see how to construct a finitely-generated
group that contains copies of these $\ell_\infty$-cubes with long side-lengths.

Finally, we remark that the methods above give a group with poor compression exponent specifically for embeddings into $L_p$ for $p < \infty$, because these are the Banach spaces to which Khot and Naor's analysis of cube-quotients in~\cite{KhoNao06} applies (see Remark 3.1 of their paper).  It is natural to ask whether some more `purely geometric' feature of the choice of Banach target could be responsible for this poor embeddability.

\textbf{Question}\quad Do the groups $(\bbZ_2^{\oplus H}/V)\ltimes H$ admit Lipschitz embeddings with positive compression exponents into any Banach space with finite cotype?

\bibliographystyle{abbrv}
\bibliography{bibfile}

\parskip 0pt

\vspace{7pt}

\small{\textsc{Department of Mathematics, Brown University, Box 1917, 151 Thayer Street, Providence, RI 09212, USA}

\vspace{7pt}

Email: \verb|timaustin@math.brown.edu|

URL: \verb|http://www.math.brown.edu/~timaustin|}

\end{document}